\documentclass[11pt,draft]{article}
\usepackage[cp1250]{inputenc}
\usepackage{amssymb}
\usepackage{amsmath}
\usepackage{fancyhdr}
\usepackage{amsthm}

\title{Radicals and  embeddings of Moufang loops in alternative loop
algebras}
\date{}
\author{Sandu N. I.}

\begin{document}
 \maketitle

\begin{abstract}
The paper defines  the notion of alternative loop algebra $F[Q]$
for any nonassociative Moufang loop $Q$ as being any non-zero
homomorphic image of the loop algebra $FQ$ of a loop $Q$ over a
field $F$. For the class $\mathcal{M}$ of all nonassociative
alternative loop algebras $F[Q]$ and for the class $\mathcal{L}$
of all nonassociative Moufang loops $Q$ are defined the radicals
$\mathcal{R}$ and $\mathcal{S}$, respectively. Moreover, for
classes $\mathcal{M}$, $\mathcal{L}$ is proved an analogue  of
Wedderburn Theorem for finite dimensional associative algebras. It
is also proved that any Moufang loop $Q$ from the radical class
$\mathcal{R}$ can be embedded into the loop of invertible elements
$U(F[Q])$ of alternative loop algebra $F[Q]$. The remaining loops
in the class of all nonassociative Moufang loops $\mathcal{L}$
cannot be embedded into  loops of invertible elements of any
unital alternative algebras.
 \maketitle
\smallskip\\
\textbf{Mathematics of subject classification}: 20N05.
\end{abstract}

\textbf{0 \ \ Introduction}

\vspace{3mm}

Embedding of a Moufang loop into a loop of invertible elements
$U(A)$ of an alternative algebra with unit $A$ (see, for example,
\cite{G}, \cite{Shest04} is one of the major questions in the
Moufang loop theory. In general, the answer to this question is
negative \cite{Shest04}, \cite{San08}. Nevertheless, many authors
look into such Moufang loops assuming that they can be embedded
into a loop of type $U(A)$ (see, Section 6). The question on
embedding a Moufang loop into a loop of type $U(A)$ is fully
solved in this paper.

To solve this question (section 6) the notion of alternative loop
algebra $F[Q]$ for any Moufang loop $Q$ is introduced. The algebra
$F[Q]$ is alternative and it is a non-zero homomorphic image of
the loop algebra $FQ$ for the loop $Q$ over a field $F$ (Section
2). Moreover, the radicals $\mathcal{R}$ and $\mathcal{S}$ are
introduced for the class $\mathcal{M}$ of all alternative loop
algebras $F[Q]$ and for the class $\mathcal{L}$ of all Moufang
loops (section 4). It is also introduced the class $\mathcal{R}_A$
into the alternative loop algebras. Let $A \in \mathcal{M}$ and
let $A = R^{\sharp} = R \oplus Fe$, i.e $A$ is obtained by
adjoining the exterior unit element $e$ to $R$. Then $A \in
\mathcal{R}_A$ when and only when $R \in \mathcal{R}$. Algebras of
type $A = R^{\sharp}$ are considered in Section 3.

The semisimple classes $\mathcal{P}$ and $\mathcal{S}$
corresponding to radicals $\mathcal{R}$ and $\mathcal{S}$
respectively are considered in section 5.

Proposition 6.3 and  Theorem 6.4  are the crucial structure
results for the examination of Moufang loops. These statements are
similar with Wedderburn Theorem for associative algebras, which is
regarded as the beginning of radical theory. In 1908 he proved
that every finite dimensional associative algebra is an extension
of the direct sum of full matrix algebras over corps with the help
of nilpotent algebra.

\textbf{Proposition 6.3.} \textit{Let $F[Q]$ be an alternative
loop algebra from the class $\mathcal{M}$ and let
$\mathcal{R}(F[Q])$ be its radical. Then  algebra
$(\mathcal{R}(F[Q]))^{\sharp} = F[G]$, $G \subseteq Q$, $F[G] \in
\mathcal{R}_A$, is nonassociative antisimple with respect to
nonassociativity or, equivalently, it does not contain subalgebras
that are nonassociative simple algebras and the quotient-algebra
$F[Q]/\mathcal{R}(F[Q])$ is a direct sum of Cayley-Dickson
algebras over their centre.}\vspace*{0.1cm}

\textbf{Theorem 6.4.} \textit{Let $\Delta$ be a prime field, let
$P$ be its algebraic closure, and let $F$ be a Galois extension
over $\Delta$ in $P$. Then the radical $\mathcal{S}(Q)$ of a
Moufang loop $Q$ is nonassociative antisimple with respect to
nonassociativity or, equivalently, it does not contain subloops
that are nonassociative simple loops and quotient-loop
$Q/\mathcal{S}(Q)$ is isomorphic to a direct product of matrix
Paige loops $M(F)$.}\vspace*{0.1cm}

For the basic properties of Moufang loops see \cite{Bruck},
\cite{CPS}, and of alternative algebras see \cite{ZSSS}.

The Cayley-Dickson algebras (simple alternative algebras) and the
Paige loops (simple Moufang loops) are quite well explored, see
\cite{ZSSS}, \cite{San11}.

Let $Q \in \mathcal{L}$. By definition $Q \in \mathcal{S}$ if and
only if $F[Q] \in \mathcal{R}_A$. Then from Proposition 6.3 and
Theorem 6.4 it follows that the construction of nonassociative
Moufang loops $Q$ by module of nonassociative simple Moufang loops
is limited to the examination of alternative loop algebras from
radical class $\mathcal{R}_A$. Such algebras are described  in
Propositions 3.3, 3.4, Corollary 3.5.

According to Lemma 6.1 and Proposition 6.2 the following
statements are equivalent for a nonassociative Moufang loop $Q \in
\mathcal{L}$:

r1) $Q \in \mathcal{S}$;

r2) the loop $Q$ is antisimple with respect to nonassociativity;

r3) the loop $Q$ does not have subloops that are simple loops.
Then the contrary statements hold:

nr1) $G \notin \mathcal{S}$, i.e. $G \in \mathcal{L} \backslash
\mathcal{S}$;

nr2) $G$ are not antisimple with respect to nonassociativity loop;

nr3) the loop $G$ contains  subloops that are simple loops hold
for any nonassociative Moufang loop $G \in \mathcal{L} \backslash
\mathcal{S}$.

From the definition of class of alternative loop algebras
$\mathcal{R}_A$, the definition of class of loops $\mathcal{S}$
and Theorem 3.2 it follows that if a nonassociative Moufang loop
$Q$ satisfies the condition r1) then the loop $Q$ can be embedded
into the loop of invertible elements $U(F[Q])$ of alternative loop
algebra $F[Q]$. On the other hand, in \cite{San08}, it was proved
that if a nonassociative  Moufang loop $G$ satisfies the condition
nr2) then the loop $Q$ is not imbedded into the loop of invertible
elements $\mathcal{U}(A)$ for a suitable unital alternative
$F$-algebra $A$, where $F$ is an associative commutative ring with
unit. As $\mathcal{S} \bigcap (\mathcal{L} \backslash \mathcal{S})
= \varnothing$ then the main result of this paper follows from the
above-mentioned statements.\vspace*{0.1cm}

\textbf{Theorem 6.5.} \textit{Any nonassociative Moufang loop $Q$
that satisfies one of the equivalent  conditions r1) - r3) can be
embedded into a loop of invertible elements $U(F[Q])$ of
alternative loop algebra $F[Q]$. The remaining loops in the class
of all nonassociative Moufang loops $\mathcal{L}$, i.e. the loops
$G \in \mathcal{L}$ that satisfy one of the equivalent conditions
nr1 - nr3 cannot be embedded into  loops of invertible elements of
any unital alternative algebras.}\vspace*{0.1cm}

From Corollary 5.10 and Theorem 6.5 it follows.\vspace*{0.1cm}

\textbf{Corollary 6.6.} \textit{Any commutative Moufang loop $Q$
can be embedded into a loop of invertible elements $U(F[Q])$ of
alternative loop algebra $F[Q]$.}\vspace*{0.1cm}

Recently a series of papers have been published, which look into
the  Moufang loops with the help of the powerful instrument of
group theory, in particular finite group theory (see, for example,
\cite{Doro}, \cite{Glaub68}, \cite{Glaub}, \cite{GrShest09},
\cite{Lieb}). For this purpose the  correspondence between Moufang
loops and groups with triality \cite{Doro} is used. The proofs
based on the correspondence are complex and cumbersome.

This paper offers another simpler approach method: to use the
Theorem 6.5, Corollary 6.6 instead of the the correspondence
between Moufang loops and groups with triality. Such examples are
presented into the end of section 6 and section 7. In section 7
this is proved on the basis of  Theorem 6.5 of the known results
from \cite{Glaub}: every finite Moufang $p$-loop is centrally
nilpotent. Paper  \cite{Glaub}  introduces the notion of group
with triality. We note that Theorem 6.5 was also used in
\cite{San99} for  proving the next statements.\vspace*{0.1cm}

If three elements $a, b, c$ of Moufang loop $Q$ are tied by the
associative law $ab\cdot c = a\cdot bc$, then they generate an
associative subloop (Moufang Theorem).\vspace*{0.1cm}

The intersection of the terms of the lower central series of a
free Moufang loop $L_X(\frak M)$ is the unit loop.\vspace*{0.1cm}

Any finitely generated  free Moufang loop is
Hopfian.\vspace*{0.1cm}

We will examine only  nonassociative Moufang loops and
nonassociative alternative algebras over a fixed field $F$.
Particularly, the alternative loop algebra $F[Q]$ corresponding to
nonassociative loop $Q$ is nonassociative. If the loop $Q$ is
commutative then algebra $F[Q]$ is also commutative. Then, in the
commutative case, we will consider that \text{char} $F = 0$ or $3$
as there are no nonassociative commutative alternative algebras
over fields of characteristic $\neq 0; 3$
\cite{ZSSS}.\vspace*{0.1cm}

Any algebra $A$ with unit $e$ is always considered nontrivial by
definition, therefore all  such algebras $A$ contain one
dimensional central subalgebra $Fe = \{\alpha a \vert \alpha \in
F\}$ (with the same unit $e \neq 0$, which allows to identify $Fe$
and $F$).

If $J$ is an ideal of algebra $A$ and the quotient algebra is an
algebra with unit $e$, then $J$ is a proper ideal of $A$ ($J \neq
A$) and $e \notin J$.  Besides, by definition, the homomorphisms
of algebras with unit $e$ is always \textit{unital}, i.e. keep the
unit. Hence, if $\varphi: A \rightarrow B$ is a homomorphism of
algebras with unit $e$, then $\ker \varphi$ is a proper ideal of
$A$, as $\varphi e = e \neq 0$.\vspace*{0.1cm}

Let $A$ be an associative algebra. By $I(A)$ we denote the set of
such elements $u \in A$ that $$u + v + uv = 0, \quad u + v + vu =
0 \eqno{(1)}$$ for some $v \in A$ and by $J(A)$ we denote the set
of all \textit{quasiregular elements} of $A$, i.e. the set of such
elements $a \in A$ that $$a + b - ab = 0, \quad a + b - ba = 0
\eqno{(2)}$$ for some $b \in A$. In the past, almost
simultaneously, with various goals on the elements of the set
$I(A)$ different authors (see., for example, \cite{Fost},
\cite{Malc}, \cite{Perl}) have introduced  the group operation
$(\otimes)$:
$$u\otimes v = uv + u + v.\eqno{(3)}$$ However, at present the so-called
circle operation $(\circ)$:
$$a\circ b = a + b - ab \eqno{(4)}$$ on the elements of the set $J(Q)$,
as in such a case the strong instrument of the theory of
quasiregular associative algebras (see., for example,
\cite{Jacob}, \cite{Ryab}) can be used to  consider operation
$(\circ)$ (we mentioned that this  paper is influenced by
\cite{Ryab}).

The operations $(\otimes)$ and $(\circ)$ on the sets $I(A)$ and
$J(A)$ respectively are groups. Thus, this paper established by
analogy a link between alternative algebras and Moufang loops. In
\cite{Good87} this link is established with the help of the
operation, defined by (3), but we believe it is not successful. In
such a form it is impossible to use the developed theory of
quasiregular alternative algebras, though in \cite{Good87} the
elements defined by (1) are wrongly called quasiregular. According
to \cite{ZSSS} a quasiregular alternative algebra can be
characterized (defined) as alternative algebra $A$ satisfying the
property that the set $A$ form a loop with respect to circle
operation $(\circ)$, defined by (4). Then the set $J(A) \subseteq
A$ define the Zhevlakov (quasiregular) radical of alternative
algebra $A$, the analog of Jacobson radical of associative algebra
theory. In \cite{Good87} the notion of Zhevlakov radical is
defined with respect to operation $(\otimes)$ that not correspond
\cite{ZSSS}. However, the loops $(I(A), \otimes)$ and $(J(A),
\circ)$ are isomorphic  by Corollary 1.7.

Let $A$ be an alternative algebra. Unlike \cite{Good87}, this
paper establishes a link between alternative algebras and Moufang
loops with the help of relation (4). In such a case the developed
theory of quasiregular alternative algebras can be used. For
example, in \cite[Theorem 1]{Good87}, it is quite cumbersomely
proved that groupoid $(I(A),\otimes)$ is a Moufang loop, but this
paper proves such a result for $(J(A),\circ)$ (Proposition 1.4)
quite easily. We will call loop $(J(A), \circ)$ \textit{circle
loop} of algebra $J(A)$. The paper also gives a full answer to the
modified question from \cite{Good87} about embedding a Moufang
loop into circle loop of a suitable alternative algebra
(Corollaries 6.8, 6.9).

\section{Circle Moufang loops}

Let $A$ be an algebra over a field  $F$. Let's consider that the
field $F$ is a module over itself. The unit $e$ of $F$ is the
generating element of the $F$-module $Fe$. We consider the direct
sum $A^{\sharp} = A \oplus Fe$ of the modules  $A$ and $Fe$ and
define on it the multiplication:
$$(a + \alpha\cdot e)(b + \beta\cdot e) = (ab + \alpha b + \beta
a) + \alpha\beta\cdot e$$ where $a, b \in A,\ \alpha, \beta \in
F$. It is easy to see that $e$ is the unit of algebra $A^{\sharp}$
and $A$ is an ideal of $A^{\sharp}$. $A^{\sharp}$ is called
\textit{the algebra obtained by adjoining the exterior unit
element $e$ to $A$}.

Consequently, it is always possible to pass, from any algebra $R$
to algebra $R^{\sharp} = R \oplus Fe$ with externally attached
unit $e$ and $R $ be an ideal of algebra $R^{\sharp} $. In
general, it is not always possible to restore algebra $R$ from the
algebra $R^{\sharp}$: it is possible that $R^{\sharp} _1 = A =
R^{\sharp}_2$, though the algebras $R_1$ and $R_2$ are not
isomorphic. However, if the algebras are given $R$ and $A =
R^{\sharp}$ then for any algebra $B$ with unit $\epsilon$ every
homomorphism $\varphi: R \rightarrow B$ unequivocally proceeds up
to homomorphism $\varphi: A \rightarrow B$ by rule:
$\varphi(\alpha e + r) = \alpha \epsilon + \varphi r$.
Particularly, the homomorphism $\pi: A \rightarrow F$, defined by
$\pi(\alpha e + r) = \alpha$, will be the only unital homomorphism
of algebra $A = R^{\sharp}$ in algebra $F \equiv F e$, continuing
the null homomorphism  of algebra $R$. Moreover,
hold.\vspace*{0.1cm}

\textbf{Lemma 1.1.} \textit{An algebra  $A$ with unit  $e$ will be
an algebra with  externally adjoined  unit (i.e. $A = R^{\sharp}$
for some algebra  $R$) when and only when there exists a
homomorphism $\pi: A \rightarrow \pi A = F$ of algebra $A$. In
such a case $A = R^{\sharp} = R \oplus Fe$, where $R = \ker
\pi$.}\vspace*{0.1cm}

\textbf{Proof.} For the homomorphism $\pi$ with $\ker \pi = R$ we
have $A/R \cong F$. Besides, $\pi$ is identical on  $Fe \equiv F$
and according to decomposition $a = \pi a + (a - \pi a) = \alpha e
+ r$ of elements $a \in A$ and equality $\pi A = Fe$ we have $A =
R^{\sharp}$. On the other hand, $R \bigcap Fe = 0$ for any proper
ideal $R$  and hence if  $A = R^{\sharp} = R \oplus Fe$ then $A
\rightarrow A/R \cong F$ will be the only unital homomorphism
$R^{\sharp} \rightarrow F$, continuing  the null homomorphism of
algebra $R$. This completes the proof of Lemma 1.1.\vspace*{0.1cm}

An alternative algebra is an algebra in which $x\cdot xy = x^2y$
and $yx\cdot x = yx^2$ are identities. Any alternative algebra
satisfies \textit{the Moufang identity}
$$(x\cdot yx)z = x(y\cdot xz). \eqno{(5)}$$The loop, satisfying
the identity (5), is called \textit{Moufang loop}.

Let $A$ be an alternative algebra with unit $e$. The element $a
\in A$ is said to have an \textit{inverse}, if there exists an
element $a^{-1} \in A$ such that $aa^{-1} = a^{-1}a = e$. It is
well known that for an alternative algebra $A$ with the unit e the
set $U(A)$ of all invertible elements of $A$ forms a Moufang loop
with respect to multiplication \cite{Mouf}. \vspace*{0.1cm}

\textbf{Lemma 1.2} \cite{ZSSS}. \textit{Let A be an alternative
algebra. Then, the following statements are equivalent:}

\textit{a) the elements $a$ and $b$ are invertible;}

\textit{b) the elements $ab$ and $ba$ are invertible.}
\vspace*{0.1cm}

The element $a$ of alternative algebra $A$ is called
\textit{quasiregular} if it satisfies the relation (2). The
element $b$ of (2) is called \textit{quasiinverse} of $a$. An
alternative algebra is called \textit{quasiregular} if any of its
elements is quasiregular. \vspace*{0.1cm}

\textbf{Lemma 1.3} \cite{ZSSS}. \textit{The following statements
are equivalent:}

\textit{a) the element $a$ of the alternative algebra $A$ is
quasiregular with quasiinverse $b$;}

\textit{b) the element $e - a$ of the algebra $A^{\sharp}$ is
inverse with the inverse element $e - b$.} \vspace*{0.1cm}

Quasiregularity is a fundamental concept in algebra theory because
it allows to define one of the most important radicals. An ideal
is called \textit{quasiregular} if it consists entirely of
quasiregular elements. Every alternative algebra $A$ has the
largest quasiregular ideal $J(A)$ such that $A/J(A)$ has no
non-zero quasiregular ideals. This ideal $J(A)$ is called the
\textit{Zhevlakov radical} and it is, of course, like the Jacobson
radical of associative algebra theory \cite{ZSSS}.\vspace*{0.1cm}

\textbf{Proposition 1.4.} \textit{ Let $A$ be an alternative
algebra and let $J(A)$ be its Zhevlakov radical. Then the set
$J(A)$ forms a Moufang loop with respect to  operation $x \circ y
= x + y - xy$.} \vspace*{0.1cm}

\textbf{Proof.} We suppose that $x$ and $y$ are quasiregular
elements with quasiinverses $a$ and $b$ respectively. We denote $u
= e - x, v = e - y$ where $e$ is the unit of algebra $A^{\sharp}$.
From Lemma 1.3 it follows that $u^{-1} = e - a,\ v^{-1} = e - b$,
where $uu^{-1} = e,\ vv^{-1} = e$, and from Lemma 1.2 it follows
that $(uv)(v^{-1}u^{-1}) = e$. From here we get that $((e - x)(e -
y))((e - b)(e - a)) = e$, $(e - x - y + xy)(e - b - a + ba) = e$,
$(e - x\circ y)(e - b\circ a) = e$. Hence the element $e - x\circ
y$ is inverse with the element $e - b\circ a$ and from Lemma 1.3
it follows that $x\circ y$ is quasiregular with quasiinverse
$b\circ a$. Consequently, the set $J(A)$ is closed under the
operation $(\circ)$.

It is easy to see that the $0$ element of $A$ is an unit for
$(\circ)$. To prove that the set $J(A)$ forms a loop under
$(\circ)$, it sufficient to show that $(x\circ y)\circ b = x$ and
similarly that $a\circ(x\circ y) = y$. Indeed, according to Lemma
1.3 $y = e - b$. Then by associativity of alternative algebras
(\cite{ZSSS}) we get $x\cdot yb = x((e - b)b) = xb - xb \cdot b =
xb - x \cdot bb = x(b - bb) = x((e -b)b) = x \cdot yb$, i.e. $x
\cdot yb = xy \cdot b$. Further, by (2) $-yb + y + b = 0$, then
$(x\circ y)\circ b = ( -xy + x + y)\circ b = xy\cdot b - xb - yb -
xy + x + y + b = xy\cdot b - xb - xy + x = x\cdot yb - xb - xy + x
= x (yb - b - y) + x = x\cdot 0 + x = x$. In this manner it is
proved that $a\circ(x\circ y) = y$. Hence $(J(A),\circ)$ is a
loop.

Finally, in order to prove the validity of Moufang identity (5) in
the loop $(J(A),\circ)$ it is sufficient to evaluate the
difference $((x\circ y)\circ x)\circ z) - x\circ(y\circ(x\circ
z))$ by (4) and by using the identity (5) for the algebra $A$,
diassociativity of $A$ and the identity $xy\cdot z + yx\cdot z =
x\cdot yz + y\cdot xz$ obtained through linearization of algebra
identity $xx\cdot z = x\cdot xz$. As a result we have obtained
that this difference is $0$. This completes the proof of
Proposition 1.4. \vspace*{0.1cm}

\textbf{Corollary 1.5.} \textit{Let $J(A)$ be the Zhevlakov
radical of the alternative algebra $A$. Then the Moufang loop
$(J(A), \circ)$ is isomorphic to $(I(A), \otimes)$.}
\vspace*{0.1cm}

\textbf{Proof.} Let $x \in J(A)$. $J(A)$ is a module, then $-x \in
J(A)$. Let $-y$ be a quasiinverse for $-x$. By (2) we have $(-x) +
(-y) - (-x)(-y) = 0, -x - y - xy = 0, -(x + y + xy) = 0, x + y +
xy = 0$. Hence by (1) $x \in I(A)$, i.e. $J(A) \subseteq I(A)$.
Inversely, let $x \in I(A)$ and let $x + y + xy = 0$. Then $(-x) +
(-y) - (-x)(-y) = 0$, and by (1) $-x \in J(A)$, $x \in J(A)$, i.e.
$I(A) \subseteq J(A)$. Hence $J(A) = I(A)$.

Now for $x \in J(A)$ we define $\varphi(x) = -x$. $\varphi$ is a
one-to-one map of $J(A)$ onto $I(A)$. Moreover,
$$\varphi(x \circ y) = -(x + y - xy) = -x + (-y) + (-x)(-y) =
\varphi(x) \otimes \varphi(y),$$ so $\varphi$ is an isomorphism of
$(J(A), \circ)$ onto $(I(A), \otimes)$, as required. $\Box$

\vspace{3mm}

Further the Moufang loop $(J(A),\circ)$, considered in Proposition
1.4 will be called the \textit{circle loop of algebra $A$} and
denoted  by $U^{\ast}(A)$. If $A$ contains the unit $e$ then the
correspondence $e - a \rightarrow a$ maps the multiplicative loop
of simple inverse elements of $A$ upon $U^{\ast}(A)$ and, in this
case, the circle operation does not offer anything new. Therefore,
we will assume further that  algebra $A$ is without unit element.

Let now $A$ be an arbitrary alternative algebra  with externally
adjoined  unit  $e$, i.e. by Lemma 1.1 with clearly distinguished
 one-dimensional subalgebra $Fe \equiv F$ with same unit $e$. We
define the mapping $\eta: A \rightarrow A$ by the rule $\eta a = e
- a \in A$ for any $a \in A$ (particularly, $\eta e = 0 \neq e =
\eta 0)$. Then, from the definition of circle loop (Proposition
1.4), we get the equality
$$(e - a)(e - b) = e - a \circ b, \eqno{(6)}$$
which, by replacing of type $x \rightarrow \eta x$, it is
rewritten as follows:
$$(e - a) \circ (e - b) = e - ab.\eqno{(7)}$$
Obviously, $\eta = \eta^{-1}$. By replacing of type $c \rightarrow
e - c$ in Lemma 1.3, we get that an element $a \in A$ is
invertible if and only if the element $e - a \in A$ is
quasiregular. Then from (6), (7) and Proposition 1.4 it follows
that $\eta$ is an isomorphism, which connects  the group of
invertible elements $U(A)$ of algebra $A$ with  circle group of
quasiregular elements $U^{\ast}(A)$ by rule

$$U^{\ast}(A) = \eta(U(A)) = \{a \in A \vert e - a \in U(A)\}.
\eqno{(8)}$$ Hence the rule
$$a \in U^{\ast}(A) \rightarrow e - (e - a)^{-1} = -a(e - a)^{-1}
= a^{\ast} \in U^{\ast}(A) \eqno{(9)}$$ defines on $A$, by
isomorphism $\eta$, the operation $a \rightarrow a^{\ast}$
\textit{of taking the quasiinverse}, defined on $U^{\ast}(A)$ and
corresponding precisely to the operation of taking the inverse,
defined on $U(A)$, according to equality $a \circ a^{\ast} = 0 =
a^{\ast} \circ a$ and the isomorphism
$$U^{\ast}(A) \cong U(A). \eqno{(10)}$$

By definition, an alternative algebra is quasiregular if any of
its  elements is quasiregular. Consequently, from definitions of
Zhevlakov radical $J(A)$ and circle loop $U^{\ast}(A)$ it follows
that for any alternative algebra $A$ the Zhevlakov radical $J(A)$
is a loop with respect to circle operation $(\circ)$ and
$$J(A) = U^{\ast}(A).\eqno{(11)}$$
In particular, an alternative algebra $A$ is quasiregular if and
only if the algebra $A$ coincides with its circle loop
$U^{\ast}(A)$. According to (9), it means that, on algebra $A$,
there also exists the unique operation $x \rightsquigarrow
x^{\ast}$, of taking the quasiinverse, related with the basic
operations of identity $x + x^{\ast} = xx^{\ast} = x^{\ast}x$
(i.e. $r \circ r^{\ast} = 0 = r^{\ast} \circ r$ for all $r \in A$
by the construction of loop $U^{\ast}(A)$).

Hence the class of all quasiregular algebras $\mathcal{K}^{\ast}$
form a variety,  if considered with an additional unitary
operation $x \rightsquigarrow x^{\ast}$ of taking the
quasiinverse. Then, by Birkhoff Theorem, the class
$\mathcal{K}^{\ast}$ is closed with respect to the taking of
quasiregular subalgebras, of direct product of quasiregular
subalgebras  and of homomorphic images of homomorphisms of
quasiregular subalgebras. But  the class $\mathcal{K}^{\ast}$ is
also closed  in respect to usual homomorphisms, i.e. to
homomorphisms of algebras. Indeed, the following result
holds.

\vspace{3mm}

\textbf{Lemma 1.6.} \textit{Let $\varphi$ be a homomorphism of
algebra $A \in \mathcal{K}^{\ast}$. Then $\varphi A \in
\mathcal{K}^{\ast}$.}\vspace*{0.1cm}

\textbf{Proof.} Let $a \in A$. From definition of quasiregular
elements it follows that $\varphi a$ is a quasiregular element.
Let $a^{\ast}$, $(\varphi a)^{\ast}$ be their quasiinverse
elements. The homomorphism $\varphi$ is unital. Hence from (9) it
follows that $(\varphi a)^{\ast} = \varphi a^{\ast}$. Then the
homomorphism $\varphi$  saves the identity $x + x^{\ast} =
xx^{\ast} = x^{\ast}x$, distinguishing quasiregular algebras.
Consequently, $\varphi A \in \mathcal{K}^{\ast}$, as
required.\vspace*{0.1cm}

From equalities $ab + (a \circ b) = a + b$  it follows easily that
a subspace $R$ of algebra $A$ is it subalgebra then and only then
$R$ is a subgroupoid of the groupoid $(A, \circ)$. Hence the
\textit{circle subgroupoid} $(R, \circ)$ of circle groupoid $(A,
\circ)$, isomorphic to  multiplicative subgroupoid $e - R$ of
algebra $A$ is linked  with any subalgebra $R$ of alternative
algebra $A$ by rule
$$e - R = \{e -r \vert r \in R\} = \eta R. \eqno({12)}$$
Then, according with (6) - (11), with
$$U^{\ast}(R) = \eta(U(e - R)) \eqno{(13)}$$ and with the definition of
multiplicative loop of algebra, we get the circle loop
$U^{\ast}(R)$ of algebra $R$, which is a subloop of loop
$U^{\ast}(A)$ and it is isomorphic to multiplicative loop $U(e -
R)$ of loop $U(A) \cong U^{\ast}(A)$ by  (10). In particular, as
$F \equiv Fe$ then the multiplicative group $U(F) = F \backslash
\{0\}$ of field $F$ and circle group $U^{\ast}(F) = F \backslash
\{0\}$ are central subloops of loops $U(A)$ and $U^{\ast}(A)$,
respectively.

From equalities $a + x + ax = x + a + xa = 0$  it follows that the
set of all elements of some subalgebra $R$ of algebra $A$, that
have quasiinverses in $R$, is a subloop of loop $U^{\ast}(A)$,
$U^{\ast}(R) \subseteq U^{\ast}(A)$ . By (11) an alternative
algebra $A$ is quasiregular if $A$ coincides with circle loop
$U^{\ast}(A)$, $J(A) = U^{\ast}(A)$. The Zhevlakov radical $J(A)$
of any alternative algebra $A$ is hereditary, i.e. $J(R) = R
\bigcap J(A)$ for any ideal $R$ of $A$. Hence
$$U^{\ast}(R) = R \bigcap U^{\ast}(A) = \{r \in R \vert e - r \in
U(A)\}. \eqno{(14)}$$  The following result holds, too.
\vspace*{0.1cm}

\textbf{Proposition 1.7.} \textit{Let A be an alternative algebra
and $R$ be an ideal of $A$. Then $U^{\ast}(R)$ is a normal subloop
of $U^{\ast}(A)$ and $U^{\ast}(A)/U^{\ast}(R) \cong
U^{\ast}(A/R)$.} \vspace*{0.1cm}

\textbf{Proof.} Let $x, y \in U^{\ast}(A), u \in U^{\ast}(R)$ and
let $a, b$ will be the quasiinverses of $x, y$ respectively. In
the proof of Proposition 1.4 it is shown that the element $x\circ
y$ is quasiregular with a quasiinverse $b\circ a$. By definition
the subloop $U^{\ast}(R)$ is normal in $U^{\ast}(A)$ if $x\circ
U^{\ast}(R) = U^{\ast}(R)\circ x, x\circ(y\circ U^{\ast}(R)) =
(x\circ y)\circ U^{\ast}(R), (U^{\ast}(R)\circ x)\circ y =
U^{\ast}(R)\circ(x\circ y)$ for any $x, y \in U^{\ast}(A)$. Any
Moufang loop is an $IP$-loop, i.e. it satisfies the identities
$x^{-1}\cdot xy = y, yx\cdot x^{-1} = y$. Then to show that
$U^{\ast}(R)$ is normal in $U^{\ast}(A)$ it is sufficient to show
that $t_1 = (x\circ u)\circ a \in U^{\ast}(R)$ and $t_2 = ((u\circ
x)\circ y)\circ(b\circ a) \in U^{\ast}(R)$. From the
aforementioned we have that $t_1, t_2 \in U^{\ast}(A)$. Further,
$u \in R$, then by (4) $t_1 = xu\cdot a - xa - ua - xu + x + u + a
= xu\cdot a - ua - xu + u \in R$ since $-xa + x + a = 0$. We
similarly have $t_2 = r - xy\cdot ba + xy\cdot b + xy\cdot a +
x\cdot ba - xb + y\cdot ba - ya - xy - ba$, where $r \in R$. We
denote $\bar{x} = e - x,\ \bar{y} = e - y,\ \bar{a} = e - a,\
\bar{b} = e - b$, where $e$ is the unit of algebra $A^{\sharp}$.
Let us express $t_2$ in terms of $x, y, a, b$ over $\bar{x},
\bar{y}, \bar{a}, \bar{b}$ respectively. We get that $t_2 = r -
\bar{x}\bar{y}\cdot \bar{b}\bar{a} + \bar{x}\bar{a} +
\bar{y}\bar{b} - e$. By Lemmas 1.2 and 1.3 $\bar{x}\bar{y}\cdot
\bar{b}\bar{a} = \bar{x}\bar{a} = \bar{y}\bar{b} = e$. Hence $t_2
= r$. Then $t_1, t_2 \in U^{\ast}(R)$ and, consequently, the
subloop $U^{\ast}(R)$ is normal in $U^{\ast}(A)$. We proved that
the ideal $R$ of algebra $A$ induced the normal subloop
$U^{\ast}(R)$ of loop $U^{\ast}(A)$.

As  noted above, the homomorphic image of circle loop
$U^{\ast}(A)$ under homomorphism $A \rightarrow A/R$ is a circle
loop. Hence the quotient loop \break $U^{\ast}(A)/U^{\ast}(R)$ is
a circle loop. By (10) $U^{\ast}(A/R) = J(A/R)$. The Zhevlakov
radical $J(A/R)$ of a maximal ideal of $A/R$. Hence
$U^{\ast}(A/R)$ is a maximal subloop of multiplicative groupoid of
algebra $A/R$.

We will show   that the quotient loop $U^{\ast}(A)/U^{\ast}(R)$ is
isomorphic to the corresponding subloop of circle loop
$U^{\ast}(A/R)$. Indeed, if $x_1$ and $x_2$ belong to the same
coset of $U^{\ast}(A)$ modulo $U^{\ast}(R)$, then $x_1 = x_2\circ
r$ where $r$ is a quasiregular element of $R$. But $x_2\circ r =
-x_2r + x_2 + r$, consequently, $x_1 - x_2 = r - x_2r \in R$.
Conversely, if $x_1 - x_2 \in R$ and $a$ is a quasiinverse for
$x_2$, i.e. $a = x_2^{-1}$, then $x_1 = x_2 + r$ ($r \in R$), $x_2
+ a - ax_2 = 0,\ x_2^{-1}\circ x_1 = a + x_1 - ax_1 = a + x_2 + r
- ax_1 = ax_2 + r - ax_1 = ax_2 + r - ax_1 = r - a(x_1 - x_2) = r
- ar \in R$. Hence $x_2^{-1} \circ x_1 \in U^{\ast}(R)$.

We proved that $U^{\ast}(A)/U^{\ast}(R) \subseteq U^{\ast}(A/R)$
or, by (10), $J(A)/J(R) \subseteq J(A/R)$. If $B/J(R) = J(A/R)$,
then by \cite[Lemma 13, cap. 10]{ZSSS}  it follows that $B
\subseteq J(A)$. Consequently, $J(A)/J(R) \cong J(A/R)$ or, by
(19), \break $U^{\ast}(A)/U^{\ast}(R) \cong U^{\ast}(A/R)$, as it
was required.

\vspace{3mm}

\textbf{Corollary 1.8.} \textit{Let $R$ be an arbitrary non-zero
alternative algebra and let $A = R^{\sharp} = Fe \oplus R$. Then
the following results hold}:

(i) \emph{the circle loop $U^{\ast}(A)$ is a direct product of the
central subloop $U^{\ast}(F)$ and the normal subloop}
$U^{\ast}(R);$

(ii) \emph{the loop of invertible elements $U(A)$ is a direct
product of central subloop $U(F)$ and normal subloop} $U(e - R)$.

\vspace{3mm}

\textbf{Proof.} By Proposition 1.7 and (10) and the
above-mentioned result as well, $U^{\ast}(R)$, $U^{\ast}(F)$ are
normal subloops of $U^{\ast}(A)$ and $U(e - R)$, $U(F)$ are normal
subloops of $U(A)$. As $R$ is a proper ideal of algebra $A$ then
$$U^{\ast}(R) \bigcap U^{\ast}(F) = \{0\},\ \ U(e - R) \bigcap U(F) =
\{e\}. \eqno{(15)}$$ Besides, as $e \notin R$ then $R \bigcap U(F)
= \varnothing = R \bigcap U(A)$. Hence $$U(A) \subseteq A
\backslash R = \{\alpha e + r\ \vert \ 0 \neq \alpha \in F,\ r \in
R\}, \eqno({16)}$$
$$U(A) = \{\alpha(e - u)\ \vert \ \alpha \in U(F),\ u \in
U^{\ast}(R)\}. \eqno{(17)}$$ Now, Corollary 1.8 follows from (14),
(15) and (17). $\Box$

\vspace{3mm}According to Lemma 1.1 the peculiarity of algebras
with externally adjoined  unit  $A = R^{\sharp}$ is linked with
the homomorphisms $\pi: A \rightarrow F$ of the considered
algebras on the one-dimensional algebra $F \equiv Fe$. An algebra
$A$ with unit $e$ will be algebra with  externally adjoined  unit
(i.e. $A = R^{\sharp}$ for some algebra $R$) when and only when $A
\neq \bigcap\{J \lhd A \vert A/J \cong F\}$ or,  equally, the set
$S_F(A) = \bigcap\{J \lhd A \vert A/J \cong F\}$ is non-empty. As
indicated in the beginning of the section, it could be the case
that  $R_1, R_2 \in S_F(A)$, i.e. that $R^{\sharp}_1 = A =
R^{\sharp}_2$, though algebras $R_1$ and $R_2$ are quite
different. But for circle loops $U^{\ast} \cong U(e - R)$ from
Corollary 1.8 the following result holds.

\vspace{3mm}

\textbf{Corollary 1.9.} \textit{If $R_1, R_2 \in S_F(A)$ then
$U^{\ast}(R_1) \cong U^{\ast}(R_2)$.}

\section{Alternative loop algebras}

Let $F$ be a field (with unit $1$) and $Q$ be a Moufang loop with
unit $e$. We remind that, by its definition, the \textit{loop
algebra} $FQ \equiv F(Q)$ is a free $F$-module with the basis $\{g
\vert g \in Q\}$ and the product of the elements of this basis is
just their product in the loop $Q$. Any element $g \in Q$ is
identified with the element $1g$, and any element $\lambda \in F$
is identified with the element $\lambda e$. In particular, the
unit of algebra $FQ$ may be  considered both  as unit of field $F$
and as unit of loop $Q$. In this case, every homomorphism
$\varphi$ of algebra $FQ$ must be \textit{unital}, i.e. it has to
maintain the unit, $\varphi e = e$. Since $\varphi e = e \neq 0$
then $\ker \varphi$ is a proper ideal of $FQ$.

\vspace{3mm}Let $H$ be a normal subloop of the loop $Q$ and let
$\omega H \equiv \omega (H)$ be the ideal of the loop algebra
$FQ$, generated by the elements $e - h$ ($h \in H$). If $H = Q$,
then $\omega Q$ will be called the \textit{augmentation ideal} of
loop algebra $FQ$. In [1, Lemma 1] it is proved that

$$F(Q/H) \cong FQ/\omega H. \eqno{(18)}$$

By definition the Moufang loop $Q$ satisfies the Moufang identity
$(xy \cdot x)z = x(y \cdot xz)$. It is easy  to see that the loop
algebra $FQ$ does not always satisfy the Moufang identity if the
loop $Q$ is nonassociative. This is an equivalent to the fact that
the equalities
$$(a,b,c)+(b,a,c) = 0,\quad (a,b,c)+(a,c,b) = 0\quad \forall a, b, c \in Q,
\eqno{(19)}$$ where the notation $(a,b,c) = ab\cdot c - a\cdot bc$
means that the associator in algebra, does not always hold  in
loop algebra $FQ$. This means that algebra  $FQ$ is not
alternative. We remind that algebra $A$ is called
\textit{alternative} if the identities $(x,x,y) = (y,x,x) = 0$
hold  in it.

Let $I(Q)$ denote the ideal of algebra $FQ$, generated by all
elements of the left part of equalities (19). It follows from the
definition of  loop  algebra and diassociativity of Moufang loops
that $FQ/I(Q)$ will be an alternative algebra. Further for the
alternative algebra $FQ/I(Q)$ we use the notation $F[Q]$ and we
call them \textit{alternative loop algebra}.

In \cite{San88}, \cite{San99} it is proved that a free Moufang
loop $L$ is isomorphically embedded under homomorphism $\eta: FL
\rightarrow F[L]$ into loop of invertible elements of algebra
$F[L]$. If the  image  $L$ is identified with $L$ then the
following holds.\vspace*{0.1cm}

\textbf{Lemma 2.1.} \textit{Any free Moufang loop $L$ is a subloop
of the loop of invertible elements $U(F[L])$ of the alternative
loop algebra $F[L]$.}

\vspace{3mm}

From the definition of loop algebra $FL$ and Lemma 2.1 it
follows.

\vspace{3mm}

\textbf{Corollary 2.2.} \textit{Any element $u$ of the alternative
loop algebra $F[L]$ of any free Moufang loop $L$ is a finite sum
$u = \sum_{i=1}^k\alpha_ig_i$, where $\alpha_i \in F$, $g_i \in
L$.}\vspace*{0.1cm}

Further we will use the following statement proved in
\cite{San88}, \cite{San99}.\vspace*{0.1cm}

\textbf{Lemma 2.3.} \textit{Let $A$ be an alternative algebra and
let $Q$ be a subloop of the loop of invertible elements $U(A)$.
Then the restriction of any homomorphism of algebra $A$  upon $Q$
will be a loop homomorphism. Consequently, any ideal $J$ of $A$
induces a normal subloop $Q \bigcap (e + J)$ of loop $Q$.}

\vspace{3mm}

Let $H$ be a normal subloop of free Moufang loop $L$ with unit
$e$. We denote the ideal of algebra $F[L]$, generated by the
elements $e - h$ ($h \in H$) by $\omega [H]$. If $H = L$, then
$\omega [L]$ will be called the \textit{augmentation ideal} of the
alternative loop algebra $F[L]$.

For a Moufang loop $Q$, let $L$ be a free Moufang loop such that
the loop $Q$ has a presentation $Q = L/H$. We consider the mapping
$\bar{\mu}: FL \rightarrow FQ$ induced by homomorphism $\mu: L
\rightarrow L/H = Q$ by
$$\bar{\mu}(\sum ^{(FL)}_{g \in L}\alpha_gg) = \sum ^{(FL)}_{g \in L}
\alpha_g\mu(g) =  \sum^{(FQ)}_{a \in Q}\alpha_aa, \eqno{(20)}$$
where $a = \mu(g)$, $\alpha_g, \alpha_a \in F$ and $\sum^{(FL)}$
means sum in $F$-module $FL$, $\sum^{(FQ)}$ means sum in
$F$-module $FQ$. The mapping $\bar{\mu}$ is defined correctly
because  $FL$ is an $F$-module with basis $\{g$ $\vert$ $\forall g
\in L\}$, $FQ$ is an $F$-module with basis $\{a \vert \forall a
\in Q\}$ and $\mu$ is an epimorphism. Moreover, as $FL$ is a free
module, then $\bar{\mu}$ is an epimorphism of $F$-modules.
Further, let $x = \sum^{(FL)}_{g \in L}\alpha_gg$, $y =
\sum^{(FL)}_{h \in L}\beta_hh$. Then $\bar{\mu}(xy) =
\bar{\mu}(\sum^{(FL)}_{g, h \in L}\alpha_g\beta_h(gh)) =
\sum^{(FL)}_{g, h \in L}\alpha_g\beta_h\mu(gh) = \sum^{(FL)}_{g, h
\in L}\alpha_g\beta_h\mu(g)\mu(h) = \sum^{(FL)}_{g \in
L}\alpha_g\mu(g) \cdot \sum^{(FL)}_{h \in L}\beta_h\mu(h)  =
\bar{\mu}(x)\bar{\mu}(y)$. Consequently,  $\bar{\mu}: FL
\rightarrow FQ$ is a homomorphism of algebras, and by (18)
$$\ker \bar{\mu} = \omega H, \quad FQ = FL/\omega H.
\eqno{(21)}$$

Let $\mu_Q$ be the homomorphism of the alternative loop algebra
$F[L]$ induced by homomorphism $\bar{\mu}$ of loop algebra $FL$,
$\mu_Q(F[L]) = \break \bar{\mu}(FL)/\bar{\mu}(I(L)) =
FQ/\bar{\mu}(I(L))$. To exclude the null homomorphisms we will
consider that the induced homomorphism $\mu_Q$ is unital.

The algebra $F[L]$ is alternative and the algebra
$FQ/\bar{\mu}(I(L))$ is alternative, as well. In this case $I(Q)
\subseteq \bar{\mu}(I(L))$. Further, by definition, the ideal
$I(L)$ of the loop algebra $FL$ is generated  by the set
$\{(u,v,w) + (v,u,w), \break (u,v,w) + (u,w,v) \vert \forall u, v,
w \in L\}$. Since $\bar{\mu}((u,v,w) + (v,u,w)) = \break
(\mu(u),\mu(v),\mu(w)) +  (\mu(v),\mu(u),\mu(w))$,
$\bar{\mu}((u,v,w) + (u,w,v)) = \break (\mu(u),\mu(v),\mu(w)) +
(\mu(u),\mu(w),\mu(v))$ and $\mu(u), \mu(v), \mu(w) \in Q$, then
\break
 $\bar{\mu}(I(L) \subseteq I(Q)$ and $\bar{\mu}(I(L)) =
I(Q)$. Consequently,
$$\mu_Q(F[L]) = FQ/I(Q) = F[Q]. \eqno{(22)}$$

Further, according to (21) and homomorphism theorems it follows
\break $I(Q) =  \overline{\mu}_Q(I(L)) = (I(L) + \omega H)/\omega
H \cong \omega H/(\omega H \bigcap I(L))$, i.e.
$$I(Q) \cong \omega
H/(\omega H \bigcap I(L)). $$

 We denote by $"$ - $"$ the difference in loop algebra
$FL$, by $" \ominus "$ we denote the difference in alternative
loop algebra $F[L]$ and by  $\theta$ -  the restriction on ideal
$\omega H$ of natural homomorphism $\eta: FL \rightarrow FL/I(L) =
F[L]$. It is obvious  that $\ker \theta = \omega H \bigcap I(L)$
and $\theta(\omega H) = \omega H/(\omega H \bigcap I(L))$.

By definition, the ideal $\omega H$ is generated by set $\{e - h
\vert \forall h \in H\}$. From Lemma 2.1, it follows that $\eta(H)
= H$. Then the algebra $\theta(\omega H)$ is generated  by the set
$\{e \ominus h \vert \forall h \in H\}$. We have $\theta(\omega H)
= \eta(\omega H)$. Recall that we have above proved the equality
$\theta(\omega H) = \omega H/(\omega H \bigcap I(L))$. By the
homomorphisms theorem it results $\omega H/(\omega H \bigcap I(L))
\cong (\omega H + I(L))/I(L)$. Hence the ideal $(\omega H +
I(L))/I(L)$ of algebra $FL/I(L)$ is generated  by the set $\{e
\ominus h \vert \forall h \in H\}$. Consequently,
$$(\omega H + I(L))/I(L) = \omega [H]. \eqno{(23)}$$

Now, by (21), homomorphism theorems and (23) it follows
$$\mu_Q(F[L]) = \mu_Q(FL/I(L)) = \bar{\mu}(FL)/\bar{\mu}(I(L)) =$$
$$((FL + \omega H)/\omega H)/((I(L) + \omega H)/\omega H) \cong
(FL + \omega H)/(I(L) + \omega H) =$$ $$FL/(I(L) + \omega H) \cong
(FL/I(L))/((I(L) + \omega H)/I(L))= F[L]/\omega[H].$$ According to
(22), it results $\mu_Q(F[L]) = F[Q]$, i.e.

$$\ker \mu_Q = \omega [H]. \eqno{(24)}$$

The homomorphism of alternative loop algebras $\mu_Q: F[X]
\rightarrow F[Q]$ is induced by homomorphism  of loop algebras
$\bar{\mu}: FX \rightarrow FQ$ which  is induced, in its turn, by
the homomorphism of loops $\mu: X \rightarrow Q$. Then, from (20),
it follows that any homomorphism of loops $\mu: X \rightarrow Q$
induces a homomorphism of alternative loop algebras $\mu_Q: F[X]
\rightarrow F[Q]$, defined by

$$\mu_Q(\sum ^{(F[L])}_{g \in L}\alpha_gg) = \sum^{(F[Q])}_{a \in
Q}\alpha_a\mu(g) =  \sum^{(F[Q])}_{a \in Q}\alpha_aa,
\eqno{(25)}$$ where $a = \mu(g)$, $\alpha_g, \alpha_a \in F$ and
$\sum^{(F[L])}$ means the sum in the $F$-module $F[L]$,
$\sum^{(F[Q])}$ means the sum in the $F$-module $F[Q]$.

We remind that in order to exclude the case $ker \mu_Q = F[Q]$ we
assume that the homomorphism $\mu_Q$ is unital. Then from (24) the
following results.\vspace*{0.1cm}

\textbf{Proposition 2.4.} \textit{Let $L$ be a free Moufang loop,
let Q be a Moufang loop which has the presentation $Q = L/H$ such
that the homomorphism $\mu_Q$, induced by (25) by homomorphism
$\mu : L \rightarrow Q$, is unital. Then the alternative loop
algebra $F[Q]$ has the presentation} $F[Q] = F[L]/\omega
[H].$

\vspace{3mm}

\textbf{Corollary 2.5.} \textit{The alternative loop algebra
$F[Q]$ of a Moufang loop $Q$ is generated as an $F$-module by the
set $\{q \vert  q \in Q\}$.}

The statement follows from Proposition 2.4 and
(25).\vspace*{0.1cm}

Now we consider a homomorphism $\rho$ of the alternative loop
algebra $F[L]$. From Lemma 2.1, it follows that $F[L]$ is
generated as $F$-module by set $\{g\ \vert \ g \in L\}$. Then the
$F$-module $\rho(F[L])$ is generated by set $\{\rho(g) \vert g \in
L\}$. Hence any element $x \in \rho(F[L])$ has a form $x = \sum_{g
\in L}\alpha_g\rho(g)$.

By Lemma 2.3 $\rho$ induces a normal subloop $H$ of loop $L$. From
(25) it follows  that the homomorphism $\mu: L \rightarrow L/H =
Q$ induces a homomorphism of alternative loop algebras $\mu_Q:
F[L] \rightarrow F[Q]$, defined by $$\mu_Q(\sum ^{(F[L])}_{g \in
L}\alpha_gg) = \sum^{(F[Q])}_{a \in Q}\alpha_a\mu(g).$$ Since
$\mu(g) = gH = \rho(g)$, it follows $\eta(F[X]) = \mu_Q(F[X]) =
F[Q]$. Hence we proved the next result.

\vspace{3mm}

\textbf{Proposition 2.6.} \textit{Let $L$ be a free Moufang loop.
The homomorphic images of the form $\mu_Q(F[L]) = F[Q]$ are the
only alternative loop algebras.  The homomorphisms $\mu_Q$ are
unital and are induced by homomorphisms of loops $\mu: L
\rightarrow Q$ by rules (25).}\vspace*{0.1cm}

\textbf{Corollary 2.7.} \textit{Let $\varphi$ be an unital
homomorphism of alternative loop algebra $F[Q]$. Then the
homomorphic image $\varphi(F[Q])$ is a non-zero alternative loop
algebra.}\vspace*{0.1cm}

\textbf{Proof.} We consider the homomorphism $\mu_Q: F[L]
\rightarrow F[Q]$ from Proposition 2.6 and let $\varphi$ be a
homomorphism of alternative loop algebra $F[Q]$. The product
$\varphi\mu_Q$ is a homomorphism of alternative loop algebra
$F[L]$ on algebra $\varphi(F[Q])$. By Proposition 2.6,
$\varphi(F[Q])$ is an alternative loop algebra, as it was
required.

\vspace{3mm}

Let $L$ be a free Moufang loop with unit $e$. By Corollary 2.2 any
element $a \in F[L]$ has a form $a = \sum_{i=1}^k\alpha_iu_i$,
where $\alpha_i \in F$, $u_i \in L$. Let $H$ be  a normal subloop
of loop $L$ and let $\varphi: L \rightarrow L/H$ be the natural
homomorphism. It is easy to see that the mapping
$\overline{\varphi}: F[L] \rightarrow F[L/H]$, defined by rule

$$\overline{\varphi}(\sum_{g \in
L}\alpha_gg) = \sum_{g \in L}\alpha_g\varphi(g) =\sum_{g \in
L}\alpha_ggH \eqno{(26)}$$ is a homomorphism. Then it necessarily
follows $$F[L/H] \cong F[L]/\ker  \overline{\varphi}.
\eqno{(27)}$$ We assume that $F[L]/\ker \overline{\varphi}$ is an
algebra with externally adjoined unit. Then $\overline{\varphi}$
is an unital homomorphism, i.e. $\overline{\varphi}(e) = e \neq
0$. In such a case   $e \notin \ker
\overline{\varphi}$.\vspace*{0.1cm}

\textbf{Lemma 2.8.} \textit{Let $\overline{\varphi}$ be a
homomorphism defined in (26) and we assume that
$F[L]/\ker\overline{\varphi}$ is an algebra with externally
adjoined  unit. Then}

\textit{1) $h \in H$ if and only if $e - h \in \omega[H]$,}

\textit{2) $F[L/H] \cong F[L]/\omega[H]$,}

\textit{3) $\omega[H] = \ker\overline{\varphi}$.}\vspace*{0.1cm}

\textbf{Proof.} 1). As the mapping $\overline{\varphi}$ is
$F$-linear, then for $u \in F[L]$ and $h \in H$
$$\overline{\varphi}((e - h)u) = (\varphi e - \varphi h)\varphi u =
(e - H)(uH) = uH - uH = 0$$ and $$\omega[H] \subseteq \ker
\overline{\varphi}. \eqno{(28)}$$

If $g \notin H$ then $gH \neq H$ and $\overline{\varphi}(e - g) =
H - gH \neq (0)$. Hence $e - g \notin \ker \overline{\varphi}
\supseteq \omega[H]$ by (28), i.e. $e - g \notin \omega[H]$.

2). Let the ideal $\omega[H]$ of algebra $F[L]$ induces, by Lemma
2.3, the normal subloop $K = L \bigcap (e - \omega[H])$ of loop
$L$ and, hence, $F[L/K] \equiv F[L]/\omega[H]$. From the first
relation, we get $1 - K \subseteq \omega[L]$. By item 1) $K = H$,
hence $F[L/H] \equiv F[L]/\omega[H]$.

3). The isomorphism $\xi: F[L]/\omega[H] \rightarrow F[L]/\ker
\overline{\varphi}$ follows from (27) and item 2). For any element
$u \in F[L]$ we denote by $\overline{u}$ the image of $u$ into
$F[L]/\omega[H]$ and by $\overline{\overline{u}}$ we denote the
image of $u$ into $F[L]/\ker \overline{\varphi}$. Let $0 \neq u
\in \ker \overline{\varphi} \backslash \omega[H]$. As $a \notin
\omega[H]$ then $\overline{0} \neq \overline{u}$. Hence
$\xi(\overline{0}) \neq \xi(\overline{u})$,
$\overline{\overline{0}} \neq \overline{\overline{u}}$. But as $u
\in \ker \overline{\varphi}$ then $\overline{\overline{0}} =
\overline{\overline{u}}$ what is a contradiction. Hence $\omega[H]
= \ker \overline{\varphi}$. This completes the proof of Lemma 2.8.

\section{Alternative loop algebras   with  externally \break adjoined  unit}

Let now $\omega[L]$ be the augmentation ideal of the alternative
loop algebra $F[L]$ of the free Moufang loop $L \neq \{e\}$.
According to Corollary 2.2 any element $a \in F[L]$ has the form
$a = \sum_{i=1}^k\alpha_iu_i$, where $\alpha_i \in F$, $u_i \in
L$.
 We denote $R = \{\sum_{u \in L}\lambda_uu \vert \sum_{u \in
 L}\lambda_u = 0\}$. Obviously, $\omega[L] \subseteq R$. Conversely,
 if $r \in R$ and $r = \sum_{u \in L}\lambda_uu$, then $-r =
- \sum_{u \in L}\lambda_uu =  (\sum_{u \in L}\lambda_u)e - \sum_{u
\in L}\lambda_qq = \sum_{u \in L}\lambda_u(e - u) \in \omega L$,
i.e. $R \subseteq \omega[L]$. Hence $$\omega[L] =  \{\sum_{u \in
L}\lambda_uu \vert \sum_{u \in L}\lambda_u  = 0\}. \eqno{(29)}$$

From (29) it follows that $\omega[L] \bigcap L = \{\emptyset\}$.
Then the algebra $\omega[L]$ will be non-zero when and only when
$L \neq \{e\}$, i.e. when $F[L] \neq Fe = 0^{\sharp}$. In such
case for any $t \in L$ the equalities $Ft \bigcap \omega[L] = 0$,
$Ft + \omega[L] = F[L]$ hold and, by (29), $t - s \in \omega[L]$.
Then, the set $B_t(\omega[L]) = \{t - s\ \vert\ s\in L,\  t \neq s
\}$ generates the $F$-module $\omega[L]$ for any $t \in L$ as $L
\neq \{t\}$ and the set $\{t\} \bigcup B_t(\omega[L])$ generates
the $F$-module $F[L]$ by Corollary 2.2. In particular, the set
$\{e\} \bigcup B_e(\omega[L])$ generates the $F$-module $F[L]$ and
the set $B_e(\omega[L])$ generates the $F$-module $\omega[L]$.
Then $F[L]/\omega[L] \cong Fe$ and, by Lemma 1.1,
$$F[L] = Fe \oplus \omega[L]  = (\omega[L])^{\sharp}. \eqno{(30)}$$

If $u \in L$ then by Lemma 1.3 $e -  u$ is a quasiregular element.
The set $B_e(\omega[L])$ generates the $F$-module $\omega[L]$. By
\cite[Lemma 10.4.12]{ZSSS}, in an alternative algebra, the sum of
quasiregular elements is a quasiregular element. Hence the
augmentation ideal $\omega[L]$ is a quasiregular algebra and from
(30) it follows that $\omega[L]$  coincides with Zhevlakov radical
$J(F[L])$, $\omega[L] = J(F[L])$.

Hence we have proved the next result.

\vspace{3mm}

\textbf{Lemma 3.1.} \textit{Let $L$ be a free Moufang loop and let
$\omega[L]$ be the augmentation ideal of the alternative loop
algebra $F[L]$. Then}

\textit{1) $\omega[L]$ is generated as an ideal of the algebra
$F[L]$, as well as an $F$-module, by set $L^{0} = \{e - u \ \vert
\ \text{for all} u \in L\}$,}

\textit{2) $\omega[L]$ is a quasiregular algebra, i.e. $\omega[L]
= J(F[L])$, where $J(F[L])$ is the Zhevlakov radical of the
algebra $F[L]$.}

\vspace{3mm}

\textbf{Theorem 3.2.} \textit{Let $Q$ be a Moufang loop with unit
$e$ such that the alternative loop algebra $F[Q]$ is an algebra
with externally adjoined  unit  $e$. Then the loop $Q$ can be
embedded into the loop of invertible elements $U(F[Q])$ of the
alternative loop algebra $F[Q]$.}

\vspace{3mm}

\textbf{Proof.} Let $Q = L/H$, where $L$ is a free Moufang loop,
let $\varphi: L \rightarrow L/H$ be the natural homomorphism and
let $\overline{\varphi}: F[L] \rightarrow F[L/H] = F[L]/\ker
\overline{\varphi}$ be the homomorphism defined in (20) by
$$\overline{\varphi}(\sum_{g \in L}\alpha_gg) =  \sum_{g \in
L}\alpha_g\varphi(g) = \sum_{g \in L}\alpha_ggH.$$ By item 3) of
Lemma 2.8, $F[L]/\ker \overline{\varphi} \cong F[L]/\omega[L]$
and, by Proposition 2.4, $F[L]/\omega[L] \cong F[Q]$. Hence $F[Q]
\cong F[L]/\ker \overline{\varphi}$.

According to (30) $F[L] = Fe \oplus \omega[L]$. $Fe$ is a field,
then from $F[L]/\omega[L] \cong Fe$ it follows that $\omega[L]$ is
a maximal proper ideal of $F[L]$. Further, $\omega[H] \subseteq
\omega[L]$. Then it is easy to see that
$\overline{\varphi}(\omega[L]) = \Delta$ is a maximal proper ideal
of $F[Q] \cong F[L]/\omega[H]$. Since the homomorphism
$\overline{\varphi}$ is unital, then $\overline{\varphi}(Fe) = Fe$
and $F[Q]/\Delta \cong Fe$. By Lemma 1.1 $F[Q] = \Delta^{\sharp} =
Fe \oplus \Delta$.

We denote $B(\Delta) = \{e - q \vert e \neq q \in Q\}$, $Q^{0} =
B(\Delta) \bigcup \{0\} = \{ e - q\ \vert\ q \in Q\}$. By item 1)
of Lemma 3.1 the augmentation ideal $\omega[L]$ of algebra $F[L]$
is generated as ideal of $F[L]$ and as $F$-module by set
$B(\omega[L]) = \{e - g\ \vert\ g\in L,\ e \neq g \}$. Further,
$\overline{\varphi}B(\omega[L]) = \{\overline{\varphi}(e - g)\
\vert\ e \neq \overline{\varphi} g \in \overline{\varphi} L\} =
\{e - q\ \vert\ e \neq q \in Q\} = B(\Delta) = Q^{0} \backslash
\{0\}$. Hence the ideal $\Delta$ is generated, as an ideal of the
algebra $F[Q]$ and as well as an $F$-module, by the set $B(\Delta)
= Q^{0} \backslash \{0\}$. From $F[Q] = Fe \oplus \Delta$ it
follows that the algebra $F[Q]$ is generated by set $Q$. According
to Corollary 2.5, any element $a \in F[Q]$ have the form
$$\sum_{i = 1}^k\alpha_iq_i, \eqno{(31)},$$ where $\alpha_i \in F$, $q_i
\in Q$.

By item 2) of Lemma 3.1 the ideal $\omega[L]$ is a quasiregular
and $\omega[L] = J(F[L])$. From $\omega[H] \subseteq \omega[L]$ it
follows that $\omega[H]$ is a quasiregular ideal. Then $\omega[H]
=  J(\omega[H])$. According to (11) and Propositions 1.7, 2.4,
$J(F[Q]) = J(F[L]/\omega[H]) = J(F[L])/J(\omega[H]) =
\omega[L]/\omega[H] = \overline{\varphi}(\omega[L]) = \Delta$.
Hence $J(F[Q]) = \Delta$ and, by (11), $\Delta$ coincide with
circle loop $U^{\ast}(F[Q])$, i.e.
$$\Delta = U^{\ast}(F[Q]). \eqno{(32)}$$

We consider the mapping $\eta: F[Q] \rightarrow F[Q]$ defined by
the rule $\eta u = e - u, \forall u \in F[Q]$. From (11) and (32)
it follows that the rule $\eta_Q: \eta_Q b = e - b$, $\forall b
\in \Delta$, defines an isomorphism between the circle loop
$(\Delta, \circ)$ and the loop of invertible elements $U(F[Q])$ of
algebra $F[Q]$ because $\eta^{-1} = \eta$. Particularly, the
restriction $\overline{\eta}_Q$ of $\eta_Q$ on $Q^{0}$ is an
isomorphism of subloop $(Q^{0}, \circ) \subseteq (\Delta, \circ)$
and  loop $Q$ defined by rule: $\overline{\eta}_Q b = e - b^{0}$,
$\forall b^{0}\in Q^{0}$. Consequently, the given Moufang loop $Q$
is a subloop of multiplicative loop of invertible elements
$U(F[Q])$ of algebra $F[Q]$. This completes the proof of Theorem
3.2.\vspace*{0.1cm}

Now we define the class $\mathcal{R}_A$ of alternative
$F$-algebras. Any alternative loop algebra  with  externally
adjoined unit $F[Q]$ belong to class $\mathcal{R}_A$. Remind that
if $F[Q] \in \mathcal{R}_A$ then the loop $Q$, the field $F$ and
the algebra $F[Q]$ have the same unit $e$ and $\varphi(e) = e$ for
any homomorphism $\varphi$ of $F[Q]$. Let $F[Q] \in
\mathcal{R}_A$. Then, from Theorem 3.2, it follows that $Q
\subseteq U(F[Q])$. This fact suggest us to give the following
definition. Let $F[Q]$ be an alternative loop algebra and let $H$
be a normal subloop of loop $Q$ such that $H \subseteq U(F[Q])$.
In such a case, we denote by $\omega[H]$ the ideal of $F[Q]$
generated by the set $\{e - h \vert h \in H\}$. If $H = Q$, then
$\omega[Q]$ will be called an \textit{augmentation ideal} of the
alternative loop algebra $F[Q]$.

\vspace{3mm}

\textbf{Proposition 3.3.} \textit{Let $\omega[Q]$ be the
augmentation ideal of  an alternative loop algebra $F[Q]$ with
externally adjoined  unit  $e$, i.e. let $F[Q] \in \mathcal{R}_A$.
Then:}

\textit{1) any element $a \in F[Q]$ has the  form $\sum_{i =
1}^k\alpha_iq_i$, where $\alpha_i \in F$, $q_i \in Q$;}

\textit{2) $F[Q] = (\omega[Q])^{\sharp} = \omega[Q] \oplus Fe$;}

\textit{3) if $F[Q] = R^{\sharp}$, then $R = \omega[Q]$;}

\textit{4) $\omega[Q]$ is generated as $F$-module by the set
$B(\omega[Q]) = \{e - q \vert e \neq q \in Q\}$;}

\textit{5) any isomorphism $\varphi$ of algebra $F[Q]$ induces the
identical isomorphism on loop $Q$ and on ideal $\omega[Q]$, as
well;}

\textit{6) $\omega[Q] =  \{\sum_{q \in Q}\alpha_qq \vert \sum_{q
\in Q}\alpha_q  = 0\}$;}

\textit{7) the algebra $\omega[Q]$ is quasiregular and coincides
with the Zhevlakov radical $J(F[Q])$, $\omega[Q] = J(F[Q])$;}

\textit{8) $\omega[Q]$ coincides with the circle loop
$U^{\ast}(F[Q])$ ($\omega[Q] = U^{\ast}(F[Q])$), i.e., by (9), on
the algebra $\omega[Q]$ there exists and it is unique, the unary
operation $x \rightsquigarrow x ^ {\ast} $ of taking the
quasiinverse, that is connected with the basic operations by
identity $x + x ^{\ast} = xx^{\ast} = x^{\ast}x$ (i.e. $r \circ
r^{\ast} = 0 = r^{\ast} \circ r$ for all $r \in \omega[Q]$ from
the construction of loop $U^{\ast}(\omega[Q])$). The circle loop
$U^{\ast}(F[Q])$ is isomorphic to the loop $U(F[Q])$ of invertible
elements under the isomorphism $\eta: u \to e - u$, $\forall u \in
\omega[Q]$. The subloop $Q^0 = B(\omega[Q]) \bigcup \{0\} = \{e -
q\ \vert\ q \in Q\}$ of the loop $U^{\ast}(F[Q])$ is isomorphic to
the given  loop $Q$, i.e. $\eta(Q^0) = Q$;}

\textit{9)  $F[Q] \backslash \omega[Q] = U(F[Q])$, i.e. the
algebra $\omega[Q]$ coincides with set of all non-invertible
elements of algebra $F[Q]$;}

\textit{10) $J(I) = I \bigcap J(F[Q])$, $U^{\ast}(I) = I \bigcap
U^{\ast}(F[Q])$ for all ideals $I$ of $F[Q]$.}

\vspace{3mm}

\textbf{Proof.} The item 1) was already proved (see the proof of
Theorem 3.2 for the equality (31)).

Theorem 3.2  is proved by showing  that the ideal $\Delta$
coincides with augmentation ideal $\omega[Q]$. Then the statements
2), 4), 7), 8) are contained in proof of Theorem 3.2.

3) Let $F[Q] = R^{\sharp} = R \oplus Fe$. By items 2), 7), 8)
$F[Q] = (\omega[Q])^{\sharp} = \omega[Q] \oplus Fe$, $\omega[Q] =
J(Q)$, $(\omega[Q], \circ) = U^{\ast}(\omega[Q])$ and by Corollary
1.9 $U^{\ast}(\omega[Q]) \cong U^{\ast}(R)$. An alternative
algebra is quasiregular if it coincides with its circle loop.
Hence the ideal $R$ is  quasiregular. $Fe$ is a field. Then from
relation $F[Q]/R \cong Fe$ it follows that $R$ is a maximal ideal
of $F[Q]$. Hence $R = J(R)$. As $J(R) = J(\omega[Q])$ then $R =
J(R) = J(\omega[Q]) = \omega[Q]$, as required.

5) As $\varphi$ is an isomorphism, then $\ker \varphi = \{0\}$. By
Lemma 2.3 it follows that $\varphi$ induces the normal subloop $Q
\bigcap (e + \ker \varphi) = e$ of loop $Q$. Hence $\varphi$
induces on $Q$ the identical isomorphism $\epsilon$. By item 4)
the ideal $\omega[Q]$ is generated   as $F$-module by set
$B(\omega[Q]) = \{e - q \vert e \neq q \in Q\}$. But $\varphi(e -
q) = e - \epsilon q = e - q$. Consequently, again by item 4),
$\varphi(\omega[Q]) = \omega[Q]$.

Using the item 1) the statement 6) is proved similarly as equality
(29).

By item 7) $U^{\ast}(\omega[Q]) = \omega[Q]$. Then the item 9)
follows from  the description of the loop $U(A)$ by equalities
(16), (17).

According to (11), the item 10) is just the equality (14). This
completes the proof of Proposition 3.3.$\Box$

\vspace{3mm}

Let $H$ be a normal subloop of the Moufang loop $Q$ and let $F[Q]
\in \mathcal{R}_A$. By item 2) of Proposition 3.3, $e \notin
\omega[Q]$. Then $e \notin \omega[H]$ and hence $F[H] \in
\mathcal{R}_A$.

Let us determine the homomorphism of $F$-algebras $\varphi$: $F[Q]
\to F[Q/H]$ by the rule $\varphi(\sum \alpha_qq) = \sum
\alpha_qHq$. The following result holds.

\vspace{3mm}

\textbf{Proposition 3.4.} \textit{Let $F[Q] \in \mathcal{R}_A$ and
let $H, H_1, H_2$ be normal subloops of the loop $Q$. Then:}

\textit{1) $\ker \varphi = \omega[H]$;}

\textit{2) $e - h \in \omega[H]$ if and only if $h \in H$;}

\textit{3) if the family of elements $\{h_i\}$??? generates the
subloop $H$, then the family of elements $\{e - h_i\}$???
generates the ideal $\omega[H]$;}

\textit{4)??? if $H_1 \neq H_2$, then  $\omega[H_1] \neq
\omega[H_2$]; if $H_1 \subset H_2$, then $\omega[H_1] \subset
\omega[H_2]$; if $H = \{H_1,H_2\}$???, then $\omega[H] =
\omega[H_1] + \omega[H_2]$;}

\textit{5) $F[Q]/\omega[H] \cong F[Q/H],\quad \omega[Q]/\omega[H]
\cong \omega[Q/H]$.} \vspace*{0.1cm}

\textbf{Proof.} The statement 1) is proved similarly with item 3
of Lemma 2.8. To prove it is necessary only to use Theorem 3.2
instead of   Lemma 2.1 and, in particular, to use the item 1) of
Proposition 2 instead  of Corollary 2.2.

2). If $q \notin H$, then $Hq \neq H$. Consequently, $\varphi(e -
q) = H - Hq \neq 0$, i.e., by  1), $e - q \notin \ker \varphi =
\omega[H]$.

3), 4). Let elements $\{h_i\}$ generate subloop $H$ and let $I$ be
the ideal, generated by the elements $\{e - h_i\}$. Obviously $I
\subseteq \omega[H]$. Conversely, let $g \in H$ and let $g =
g_1g_2$, where $g_1, g_2$ are words from $\{h_i\}$. We suppose
that $e - g_1, e - g_2 \in I$. Then $e - g = (e - g_1)g_2 + e -
g_2 \in I$,  i.e. $\omega[H] \subseteq I$ and $I = \omega[H]$. Let
$H_1 \neq H_2$
 (respect. $H_1 \subset H_2$) and $g \in H_1, g \notin H_2$. Then,
 by  item 2), $e - g \in \omega[H_1]$, but $e - g \notin \omega[H_2]$.
 Hence $\omega H_1 \neq \omega
H_2$ (respect. $\omega H_1 \subset \omega H_2)$. If $H = \{H_1,
H_2\}$, then by the first statement of  3), $\omega H = \omega H_1
+ \omega H_2$.

5). Mapping $\varphi:F[Q] \to F[Q/H]$ is the homomorphism of
alternative loop algebras and, as by item 1), $\ker \varphi =
\omega[H]$, then $F[Q/\omega[H] \cong F[Q/H]$. The mapping
$\omega[Q] \to \omega[Q]/\omega[H]$ save the sum of coefficients
then from item 5) of Proposition 3.3 it follows that $\omega
L/\omega H \cong \omega (L/H)$.\vspace*{0.1cm}

\textbf{Corollary 3.5.} \textit{For a normal subloop $H$ of a
Moufang loop $Q$  with unit $e$ the following statements are
equivalent:}

\textit{1) $F[H] \in \mathcal{R}_A$;}

\textit{2) $F[H] = \omega[H] \oplus F1$;}

\textit{3) $e \notin \omega [H]$;}

\textit{4) $\omega[H]$ is a proper ideal of algebra $F[Q]$;}

\textit{5) $\omega [H] = \{\sum_{q \in H}\alpha_qq \vert  \sum_{q
\in H}\alpha_q = 0\}$.}\vspace*{0.1cm}

\textbf{Proof.} The equivalence of items 1), 2) follows from Lemma
1.1 and item 3) of Proposition 3.3. The implications 2)
$\Leftrightarrow$ 3), 3) $\Leftrightarrow$ 4) are obvious. The
implication 1) $\Leftrightarrow$ 5) follows from items 2), 6) of
Proposition 3.3.

\section{Radicals in alternative loop algebras and\\ Moufang loops}

Let $\mathcal{M}$ denote the class of all alternative loop
algebras and its augmentation ideals. If $\varphi$ is a non-zero
homomorphism of algebra $F[Q] \in \mathcal{M}$ then, by Corollary
2.7, the homomorphic image $\varphi(F[Q])$ will be an alternative
loop algebra. Hence $\varphi(F[Q]) \in \mathcal{M}$.

Let now $\omega[Q]$ be the  ideal of the alternative loop algebra
$F[Q] \in \mathcal{M}$ with unit $e$, defined above, and let
$\psi$ be a homomorphism of $\omega[Q]$. If $e \in \omega[Q]$ then
$\omega[Q] = F[Q]$ and $\psi$ is the zero homomorphism. We suppose
that $e \notin \omega[Q]$. Then $F[Q] \in \mathcal{R}$ and, by
Corollary 3.5, $F[Q] = \omega[Q] \otimes Fe$. We extend $\psi$ to
the homomorphism $\varphi$ of $F[Q]$ considering that $\varphi(e)
= e$.  Let $\varphi(F[Q]) = F[G]$. Then $F[G] = \psi(\omega[Q])
\otimes \varphi(Fe) = \psi(\omega[Q]) \otimes Fe$. By Lemma 1.1 it
follows $F[G] \in \mathcal{R}$.  Then by item 3) of Proposition
3.3 $\psi(\omega[Q] = \omega[G]$. Consequently, we proved that the
class $\mathcal{M}$ is closed with respect to homomorphic images.

Let $J$ be an ideal of the alternative loop algebra $F[Q]$. By
Lemma 2.3 $J$ induces the normal subloop $H = Q \bigcap (e + J)$
of loop $Q$. In its turn $H$ induces the homomorphism $\varphi:
F[Q] \rightarrow F[Q]/H$ defined by $\varphi(\sum\alpha_qq) =
\sum\alpha_qqH$, $\alpha_q \in F$, $q \in Q$. We have $\ker
\varphi = J$. According to item 1) of Proposition 3.4 $J = \ker
\varphi = \omega[H]$. If $J_1$ is another ideal of $F[Q]$, $J_1
\neq J$, then $\omega[H] = J \neq J_1 = \omega[H_1]$, where $H_1 =
Q \bigcap (e + J_1)$.  By item 3) of Proposition 3.4, it follows
that $H_1 \neq H$. Consequently, various ideals of algebra $F[Q]$
induce various normal subloops of the loop $Q$.

Conversely, let $H \neq H_1$ be normal subloops of loop $Q$. The
subloops $H, H_1$ induce the homomorphisms $\varphi, \varphi_1$ of
algebra $F[Q]$ and, by items 1), 3)  of Proposition 3.4, $\ker
\varphi = \omega[H] \neq \omega[H_1] = \ker \varphi_1$. Hence
various normal subloops of the loop $Q$ induce various ideals of
algebra $F[Q]$. The proper ideals of algebra $F[Q]$ have the form
$\omega[H]$, where $H$ is a normal subloop of the loop $Q$ and
$\omega[H]$ is the augmentation ideal without unit of the
alternative loop algebra $F[H]$. Consequently, the correspondence
$\omega[H] \rightarrow H$ is an one-to-one mapping between all
normal subloops $H$ of the loop $Q$ and all ideals of the algebra
$F[Q]$. Further let's consider, that all considered algebras
belong to class $\mathcal{M}$, i.e. they have the form $F[Q]$, any
of its ideals $J \neq F[Q]$ has the form $\omega[H]$, where
$\omega[H]$ is the augmentation ideal of some alternative loop
algebra $F[H]$ where $H$ is a normal subloop of loop $Q$.

Now we consider the class of algebras $\mathcal{R}_A \subseteq
\mathcal{M}$, defined in section 3. This class was analysed in
Propositions 3.3, 3.4 and Corollary 3.5. The class $\mathcal{R}_A$
is characterized  by the property that any algebra from
$\mathcal{R}_A$ is an alternative loop algebra $F[Q]$ such that
its ideals $J \neq F[Q]$ are augmentation ideals $\omega[H]$
without unit of some alternative loop algebras $F[H]$ from
$\mathcal{R}_A$, where $H \subset Q$ is a normal subloop of the
loop $Q$.

We denote by $\mathcal{R}$ the class of such augmentation ideals
$\omega[Q]$ of its alternative loop algebra $F[Q] \in
\mathcal{R}_A$.

An ideal $I$ of algebra $A \in \mathcal{M}$ will be called
\textit{$\mathcal{R}$-ideal} if it belongs to class $\mathcal{R}$.
An alternative loop algebra $F[Q] \in \mathcal{M}$, containing
non-zero ideals  will be called \textit{$\mathcal{R}_A-algebra$}
if its augmentation ideal $\omega[Q]$ belongs to class
$\mathcal{R}$, i.e. the ideal $\omega[Q]$ is without unit and
$(\omega[Q])^{\sharp} = F[Q]$.\vspace*{0.1cm}

\textbf{Lemma 4.1.} \textit{Any algebra $A$ of class $\mathcal{M}$
contains a unique maximal $\mathcal{R}$-ideal $\mathcal{R}(A)$.
The ideal $\mathcal{R}(A)$ coincides with augmentation ideal
$\omega[K]$ of some alternative loop algebra $F[K] \in
\mathcal{R}_A$.}\vspace*{0.1cm}

\textbf{Proof.} Let $\Sigma$ denote the set of all
$\mathcal{R}$-ideals of algebra $A$. The set $\Sigma$ is
non-empty, as the ideal $(0) \in \Sigma$. By Zorn Lemma the set
$\Sigma$ contains a maximal $\mathcal{R}$-ideal $\mathcal{R}(A)$.

Let's show that $\mathcal{R}(A)$ is an unique maximal
$\mathcal{R}$-ideal. Let $\mathcal{R}(B)$ also be a maximal
$\mathcal{R}$-ideal and let $x \in \mathcal{R}(A)$, $y \in
\mathcal{R}(B)$. As $\mathcal{R}(A), \mathcal{R}(B) \in
\mathcal{R}$, then, by item 5) of Corollary 3.5, $x = \sum_{g \in
Q}\alpha_gg$ with $\sum_{g \in Q}\alpha_g = 0$ and $y = \sum_{g
\in Q}\beta_gg$ with $\sum_{g \in Q}\beta_g = 0$. Then $x + y =
\sum_{g \in Q}\gamma_gg$ with $\sum_{g \in Q}\gamma_g = 0$. Hence
$\mathcal{R}(A) + \mathcal{R}(B) \in \mathcal{R}$. If
$\mathcal{R}(A) \neq \mathcal{R}(B)$ then $\mathcal{R}(A) +
\mathcal{R}(B)$ strictly contain $\mathcal{R}(A)$. Contradiction.
Consequently, $\mathcal{R}(A)$ is the unique maximal
$\mathcal{R}$-ideal.

The second statement of lemma follows from relation
$\mathcal{R}(A) \in \mathcal{R}$ and the construction of ideals of
class $\mathcal{R}$. This completes the proof of Lemma
4.1.\vspace*{0.1cm}

\textbf{Theorem 4.2.} \textit{The class $\mathcal{R}$ of all
augmentation ideals without unit is radical in class $\mathcal{M}$
of all alternative loop $F$-algebras and its
ideals.}\vspace*{0.1cm}

\textbf{Proof.} According to the definition of radical \cite{ZSSS}
should prove  the statements:

(a) any homomorphic image of any $\mathcal{R}$-ideal is an
$\mathcal{R}$-ideal;

(b) each algebra $A$ from  $\mathcal{M}$ contains an
$\mathcal{R}$-ideal $\mathcal{R}(A)$, containing all
$\mathcal{R}$-ideals of algebra $A$;

(c) the quotient-algebra $A/\mathcal{R}(A)$ does not contain any
non-null $\mathcal{R}$-ideals.\vspace*{0.1cm}

Really, let $\omega[H] \in \mathcal{R}$ and let $x = \sum_{h \in
H}\alpha_hh \in \omega [H]$. By item 5) of Corollary 3.5, $\sum_{h
\in H}\alpha_h = 0$. Any homomorphism $\varphi$ of the ideal
$\omega [H]$ does not change the sum of coefficients, $\sum
\alpha_h$. From here, it follows that $\varphi(\omega [H]) \in
\mathcal{R}$ and the statement (a) is proved.

The Lemma 4.1 is just the statement (b).

Let $A \in \mathcal{M}$. The homomorphism $\varphi: A \rightarrow
A/\mathcal{R}(A)$ maintains the sum of coefficients. Hence if $J
\neq (0)$ is an $\mathcal{R}$-ideal of $A/\mathcal{R}(A)$ then the
inverse image $\varphi^{-1}J$ will be an $\mathcal{R}$-ideal and
$\mathcal{R}(A) \subset \varphi^{-1}J$. But this contradicts the
maximality of $\mathcal{R}$-ideal $\mathcal{R}(A)$. Consequently,
$J = (0)$ and the statement (c) is proved. This completes the
proof of Theorem 4.2.\vspace*{0.1cm}

Let $A \in \mathcal{M}$. The mapping $A \rightarrow
\mathcal{R}(A)$ is called radical, defined in class of algebra
$\mathcal{M}$; denote it by $\mathcal{R}$.

Let now introduce some notions, derived from the general concepts
of the theory of radicals \cite{ZSSS}. \textit{The ideal}
$\mathcal{R}(A)$ of algebra $A \in \mathcal{M}$ is called its
\textit{$\mathcal{R}$-radical}. \textit{An alternative loop
algebra} $F[Q] \in \mathcal{M}$ will be called
\textit{$\mathcal{R}$-radical} if $F[Q] \in \mathcal{R}_A$ and
$\omega[Q] = \mathcal{R}(F[Q])$. \textit{Non-zero algebras} $F[Q]
\in \mathcal{M}$, whose radical is null, will be called
\textit{$\mathcal{R}$-semisimples}. The class $\mathcal{P}$ of all
$\mathcal{R}$-semisimple algebras of class $\mathcal{M}$ is called
\textit{semisimple class of radical} $\mathcal{R}$.

Let $A \in \mathcal{M}$. According to Lemma 4.1 the radical
$\mathcal{R}(A)$ coincides with the augmentation ideal $\omega[K]$
of some alternative loop algebra $F[K] \in \mathcal{R}_A$. Then,
by item 7) of Proposition 3.3, $\mathcal{R}(A)$ coincides with the
Zhevlakov radical $J(F[K])$ which, by item 10) of Proposition 3.3,
is hereditary. Remind that the radical $\mathcal{R}$ in the class
of algebras $\mathcal{A}$ is called hereditary if $\mathcal{R}(J)
= J \bigcap \mathcal{R}(A)$ for any algebra $A \in \mathcal{A}$
and any it ideal $J$. Then the following holds.

\vspace{3mm} \textbf{Corollary 4.3.} \textit{The radical
$\mathcal{R}$ is hereditary in class
$\mathcal{M}$.}\vspace*{0.1cm}

\textbf{Corollary 4.4.} \textit{Let $A$ be an algebra of class
$\mathcal{M}$ and let $J$ be an ideal of $A$. The following
results hold:}

\textit{(1) if $A \in \mathcal{R}$ then $J \in \mathcal{R}$;}

\textit{(2) if $A \in \mathcal{P}$ then $J \in
\mathcal{P}$.}\vspace*{0.1cm}

It follows from Corollary 1 and \cite[Theorem 3, cap.
8]{ZSSS}.

\vspace{3mm}

Let $F[Q]$ be the alternative loop algebra of a Moufang loop $Q$
with unit $e$. According to Corollary 2.5 any element $a \in F[Q]$
is a finite sum $a = \sum_{q \in Q}\alpha_qq$, where $\alpha_q \in
F$. Then we may define the ideal of $F[Q]$ generated by set $\{e -
q \vert q \in Q\}$. We denote it by $\omega[Q]$. Note that in a
similar way we have above defined the augmentation ideal
$\omega[Q]$ of algebra $F[Q] \in \mathcal{R}_A$.

\vspace{3mm}

\textbf{Corollary 4.5.} \textit{For any algebra $A$ of class
$\mathcal{P}$ the following statements hold:}

\textit{1) if $A = \omega [Q]$ then $\omega [Q] = F[Q]$;}

\textit{2) if $x \in A$ and $x = \sum_{g \in Q}\alpha_gg$ then
$\sum_{g \in Q}\alpha_g \neq 0$;}

\textit{3) any ideal $J$ of algebra $A$ has the form $J = F[H]$
and, if $J \neq 0$, then $J$ is nonassociative.}

\vspace{3mm}

\textbf{Proof.} If $A \in \mathcal{M}$, then the from definition
of class $\mathcal{M}$, it follows that $A = F[Q]$ for some
alternative loop algebra $F[Q]$. By Lemma 4.1, $\mathcal{R}(A) =
\omega[H]$. If $A \in \mathcal{P}$ then $\mathcal{R}(A) = \{0\}$.
From here it follows that the algebra $A$ does not have  non-zero
proper ideals. Then $\omega[Q] = F[Q]$.

The item 2) follows from item 6) of Proposition 3.3.

Further, $A \in \mathcal{P}$ implies $J \in \mathcal{P}$ for any
ideal $J$ of algebra $A$, by Corollary 4.4. We have above proved
that any ideal $I$ of $A$ have the form $I = \omega[H]$. Then, by
item 1), it follows $J = F[H]$. If $F[H]$ is associative then $H$
is a group and $F[H]$ is a group algebra. From the definition of
group algebra it follows that $F[H]$ is a free $F$-module with
bases $\{h \in H\}$. Then $e \notin \omega [H]$ and $F[H] \neq
\omega [H]$. Contradiction. Consequently, the ideal $J \neq 0$
cannot be associative. This completes the proof of Corollary
4.5.\vspace*{0.1cm}

In the beginning of the section we showed that for any alternative
loop algebra $F[Q]$ the mapping $\omega[H] \rightarrow H$ is an
one-to-one mapping between all normal subloops $H$ of loop $Q$ and
all ideals of algebra $F[Q]$. Moreover, the following statement
holds.\vspace*{0.1cm}

\textbf{Lemma 4.6.} \textit{Let $F[Q]$ be an alternative loop
algebra and let $H, H_1, H_2$ be normal subloops of loop $Q$.
Then:}

\textit{1) $e - h \in \omega[H]$ if and only if $h \in H$;}

\textit{2) if the elements $\{h_i\}$ generate the subloop $H$,
then the elements $\{1 - h_i\}$ generate the ideal $\omega[H]$; if
$H_1 \neq H_2$, then  $\omega[H_1] \neq \omega[H_2$]; if $H_1
\subset H_2$, then $\omega[H_1] \subset \omega[H_2]$; if $H =
\{H_1,H_2\}$, then $\omega[H] = \omega[H_1] +
\omega[H_2]$.}\vspace*{0.1cm}

\textbf{Proof.} If $e \notin \omega[H]$ then $F[H] \in
\mathcal{R}_A$ and the statement 1) is the statement 2) of
Proposition 3.4. If $e \in \omega[H]$ then $\omega[H] = F[H]$ and
statement 1) follows, from property (31), that $F[H]$ is generated
as an $F$-module by the set $H$.

The statement 2) is proved similarly as item 3) of Proposition
3.4. We have to use the item 1), only.

By $\mathcal{L}$ denote the class of all Moufang loops and by
$\mathcal{S}$ denote the class of Moufang loops $G$ such that
$F[G] \in \mathcal{R}_A$ or, equivalently, $e \notin \omega[G]$
(by Proposition 3.3). Any loop from the class $\mathcal{L}$
(respect. $\mathcal{S}$) will be called
\textit{$\mathcal{L}$-loop} (respect.
\textit{$\mathcal{S}$-loop}).  Now, let $Q \in \mathcal{L}$ be a
Moufang loop, $F[Q] \in \mathcal{M}$ be it alternative loop
algebra and, according to Theorem 4.2, let $\mathcal{R}(F[Q])$ be
the $\mathcal{R}$-radical of $F[Q]$. By Lemma 4.1,
$\mathcal{R}(F[Q]) = \omega[\mathcal{S}(Q)]$, where
$\omega[\mathcal{S}(Q)]$ is the augmentation ideal of some
alternative loop algebra $F[\mathcal{S}(Q)] \in \mathcal{R}_A$. By
Theorem 4.2 the mapping $\mathcal{R}: F[Q] \rightarrow
\mathcal{R}(F[Q]) = \omega[\mathcal{S}(Q)]$ is a radical of class
$\mathcal{M}$. Obviously, $\mathcal{R}$ induces the mapping
$\mathcal{S}: Q \rightarrow \mathcal{S}(Q)$.

Note that, with the help of Lemma 2.3, it is easy to see that from
$\mathcal{R}(F[Q]) \break = 0$ it follows $\mathcal{S}(Q) = e$.
Further, we will show that the mapping $\mathcal{S}$ is a radical
of the class $\mathcal{L}$ of loops. For this, the class of loops
$\mathcal{S}$ should satisfy the following conditions:

$\bullet$ the homomorphic image of any $\mathcal{S}$-loop is a
$\mathcal{S}$-loop;

$\bullet$ each $\mathcal{L}$-loop $Q$ contains a normal
$\mathcal{S}$-subloop $\mathcal{S}(Q)$, containing all normal
$\mathcal{S}$-subloops of the loop $Q$;

$\bullet$ the quotient loop $Q/\mathcal{S}(Q)$ does not contain
non-unitary normal \break $\mathcal{S}$-subloops.\vspace*{0.1cm}

\textbf{Theorem 4.7.} \textit{The class $\mathcal{S}$ is radical
in the class $\mathcal{L}$ of all Moufang loops.}

\vspace{3mm}

\textbf{Proof.} Let $G \in \mathcal{S}$. Then $F[G] \in
\mathcal{R}_A$. Any homomorphism $\varphi$ of loop $G$ induces a
homomorphism $\overline{\varphi}: F[G] \rightarrow F[\varphi G]$
defined by rules $\overline{\varphi}(\sum_{g \in G}\alpha_gg) =
\sum_{g \in G}\alpha_g\varphi g$, $\overline{\varphi}(F[G]) =
F[\varphi G]$. By Theorem 4.2 $\overline{\varphi}(F[G]) \in
\mathcal{R}_A$. Then $F[\varphi G] \in \mathcal{R}$. Hence
$\varphi G \in \mathcal{S}$. Consequently, the class $\mathcal{S}$
is closed under homomorphisms.

Let now $Q$ be a Moufang loop and let, by Theorem 4.2,
$\mathcal{R}(Q)$ be the $\mathcal{R}$-radical of $F[Q]$. By Lemma
2.3 the ideal $\mathcal{R}(Q)$ of $F[Q]$ induces the normal
subloop $\mathcal{S}(Q) = Q \bigcap (e + \mathcal{R}(Q))$ of the
loop $Q$. Let $\psi: F[Q] \rightarrow F[Q/\mathcal{S}(Q)]$ be the
homomorphism defined by: $\sum_{g \in Q}\alpha_gg \rightarrow
\sum_{g \in Q}\alpha_g(g\mathcal{S}(Q))$. Then, $F[Q]/\ker \psi
\cong F[Q/\mathcal{S}(Q)]$ and $\mathcal{R}(Q) \subseteq \ker
\psi$.

If $e \in \ker \psi$, then $\ker \psi = F[Q]$ and
$F[Q/\mathcal{S}(Q)] \cong F[Q]/\ker \psi = \break F[Q]/F[Q] = 0$.
We get a contradiction because $e \in F[Q/\mathcal{S}(Q)]$. Hence
$\ker \psi$ is a proper ideal of $F[Q]$. Then, as it was shown in
the beginning of section, $\ker \psi$ is an augmentation ideal.
The radical $\mathcal{R}(Q)$ is a maximal augmentation ideal of
$F[Q]$. Hence $\ker \psi = \mathcal{R}(Q)$ and
$F[Q/\mathcal{S}(Q)] \cong F[Q]/\mathcal{R}(Q)$,
$\mathcal{R}(F[Q/\mathcal{S}(Q)])\break \cong
\mathcal{R}(F[Q]/\mathcal{R}(Q))$. By Theorem 4.2,
$\mathcal{R}(Q))$ is a maximal ideal of $F[Q]$ such that the
$\mathcal{R}$-radical of the quotient-algebra
$F[Q]/\mathcal{R}(Q)$ is zero. Hence \break
$\mathcal{R}(F[Q/\mathcal{S}(Q)])  = 0$ and
$\mathcal{S}(Q/\mathcal{S}(Q)) = 0$. Consequently, the normal
sub\-loop $\mathcal{S}(Q)$ of loop $Q$ is maximal and such that
$\mathcal{S}(Q/\mathcal{S}(Q)) = 0$. This completes the proof of
Theorem 4.7.$\Box$

\vspace{3mm}Let $Q \in \mathcal{L}$. By Theorem 4.7, the mapping
$Q \rightarrow \mathcal{S}(Q)$ is a radical defined in the class
of loops $\mathcal{L}$; denote it by $\mathcal{S}$. The normal
subloop $\mathcal{S}(Q)$ of loop $Q$ will be called its
\textit{$\mathcal{S}$-radical}. \textit{A loop} coinciding with
its $\mathcal{S}$-radical will be called
\textit{$\mathcal{S}$-radical}, and \textit{the non-unitary
loops}, whose radical is equal to unit, will be called
\textit{$\mathcal{S}$-semisimples}. The class $\mathcal{T}$ of all
$\mathcal{S}$-semisimples algebras in class $\mathcal{L}$ will be
called \textit{semisimple class of radical $\mathcal{S}$}.

\vspace{3mm}

\textbf{Proposition 4.8.} \textit{The radical $\mathcal{S}$ is
hereditary in the class $\mathcal{L}$ of all Moufang loops, i.e.
for any loop $Q \in \mathcal{L}$ and its normal subloop $H$,
$\mathcal{S}(H) = H \bigcap \mathcal{S}(Q)$.}

\vspace{3mm}

\textbf{Proof.} By Theorem 4.7, the radical $\mathcal{S}(G)$ is a
maximal normal sub\-loop $H$ of the loop $G$ with respect to
property $\omega [H] \in \mathcal{R}$. From item 3) of Proposition
3.4 it follows that $\omega [\mathcal{S}(G)]$ is a  maximal ideal
$\omega [H]$ of algebra $\omega [G]$ with respect to property
$\omega [H] \in \mathcal{R}$. Then by Theorem 4.2 $\omega
[\mathcal{S}(H)] = \mathcal{R}(\omega [H])$.

Let now $H$ be a normal subloop of loop $Q$. Then $\omega [H]$
will be a normal subloop of loop $\omega [Q]$. By Corollary 4.3
$\mathcal{R}(\omega [H]) = \omega [H] \bigcap \mathcal{R}(\omega
[Q])$ and by the precious  equality $\omega [\mathcal{S}(H)] =
\omega [H] \bigcap \omega [\mathcal{S}(Q)]$. Then from  item 3) of
Proposition 4.3 it follows from here that $\omega [\mathcal{S}(H)]
= \omega [H \bigcap \mathcal{S}(Q)]$, $\mathcal{S}(H) = H \bigcap
\mathcal{S}(Q)$. This completes the proof of Proposition
4.8.\vspace*{0.1cm}

\textbf{Corollary 4.9.} \textit{Let $Q$ be a Moufang loop and let
$K$ be a normal subloop  of $Q$. Then the following statements
hold:}

\textit{(i) if $Q \in \mathcal{S}$ then $K \in \mathcal{S}$;}

\textit{(ii) if $Q \in \mathcal{T}$ then $K \in
\mathcal{T}$.}\vspace*{0.1cm}

These follow from Proposition 3 and \cite[Theorem 3, cap.
8]{ZSSS}.

\section{Semisimple alternative loop algebras and \break semisimple Moufang
loops}

Let $A$ be an algebra. \textit{The sum of ideals} $\{I_s \vert s
\in S\}$ of algebra $A$ is called the ideal $I$ of $A$ generated
by reunion $\bigcup_{s \in S}I_s$. The ideal $I$ consists of
elements $x$, presented in the form  $x = x_1 + \ldots + x_k$,
where $x_j \in I_{s_j}$ for some $s_j \in S$ and denote $I =
\sum_{s \in S}I_s$. The sum is called \textit{direct} if $I_s
\bigcap \sum_{s \neq t \in S}I_t = 0$. Denote $I =
\sum^{\oplus}I_s$ and $I = I_1 \oplus \ldots \oplus I_k$ for
finite sum of ideals.

By analogy, \textit{the product $N$ of normal subloops} $\{N_s
\vert s \in S\}$ of the loop $Q$ consists of elements $x$,
presented in the form $x = x_1 \cdot \ldots \cdot x_k$, where $x_j
\in N_{s_j}$ for some $s_j \in S$ and denote $N = \prod_{s \in
S}N_s$. The product is called \textit{direct} if $N_s \bigcap
\prod_{s \neq t \in S}N_t = 1$. Denote $N = \prod^{\otimes}N_s$
and $N = N_1 \otimes \ldots \otimes N_k$ for a finite factor
product.

\textit{An ideal} $J$ of the algebra $A$ is called \textit{simple}
if $J$ does not have other ideals  of $A$ besides the null and
\textit{ideal} $J$ itself and it is called \textit{principal} if
it is generated by one element. \textit{A normal subloop} $N$ of
loop $Q$ will be called \textit{simple} if $N$ does not have other
normal subloops of $Q$ besides the unitary subloop and loop $N$
itself.\vspace*{0.1cm}

\textbf{Lemma 5.1.} \textit{The following statements are
equivalent for  the simple normal subloops $\{N_s \vert s \in S\}$
of a loop $Q$:}

\textit{1) $Q = \prod_{s \in S}N_s$;}

\textit{2) $Q = \prod_{s \in S}^{\otimes}N_s$.}\vspace*{0.1cm}

\textbf{Lemma 5.2.} \textit{The following statements are
equivalent for the family of simple ideals $\{I_s \vert s \in S\}$
of algebra $A$:}

\textit{1) $A = \sum_{s \in S}I_s$;}

\textit{2) $A = \sum_{s \in S}^{\oplus}I_s$.}\vspace*{0.1cm}

\textbf{Proof.} The Lemmas 5.1, 5.2 are similar proofs. Let us
prove Lemma 5.2.

Let $T$ be a maximal subset of $S$ such that the sum $\sum_{t \in
T}I_t$ is direct. The sum $\sum_{t \in T}I_t$ is an ideal of $A$.
Let us show that this sum coincides with $A$. For this it is
enough to show that each ideal $I_j$ is contained in this sum. The
intersection of our sum  with $I_j$ is an ideal in  $A$ and,
consequently, equals  $0$ or $I_j$. If it equals $0$, then subset
$T$ is not maximal, as we can add $j$ to it. Consequently, $I_j$
is contained in the sum  $\sum_{t \in T}I_t$. This completes the
proof of Lemmas 5.2.\vspace*{0.1cm}

Let $A$ be an $F$-algebra and let $I(M)$ be the ideal of $A$
generated by set $M \subseteq A$. The ideal $I(M)$ consists of all
possibly types of finite sums of elements of form
$$\varphi(x_1, \ldots, x_j, a, x_{j + 1}, \ldots, x_n)_{\alpha},
\eqno{(33)}$$ where $\varphi \in F$, $a \in M$, $x_i \in A$,
$\alpha$ is a certain distribution of parenthesis.

Let $a, b \in A$. Then from (33) it follows that
$$I(a + b) \subseteq I(a) + I(b). \eqno{(34)}$$

Now, we consider an  ideal $\omega [Q]$ of the alternative loop
algebra $F[Q]$ of a Moufang loop $Q$. If $a, b \in Q$ then $e - ab
= (e - a) + (e - b) - (e - a)(e - b)$. Denote $e - u =
\overline{u}$. By (34), $I(\overline{ab}) = I(\overline{a} +
\overline{b} - \overline{a}\overline{b})$, $I(\overline{ab})
\subseteq I(\overline{a}) + I(\overline{b}) -
I(\overline{a}\overline{b})$. From (33), it follows that
$I(\overline{a}\overline{b}) \subseteq I(\overline{a})$. Then
$$I(\overline{ab}) \subseteq I(\overline{a}) + I(\overline{b}).
\eqno{(35)}$$ Moreover, the following result holds.\vspace*{0.1cm}

\textbf{Lemma 5.3.} \textit{Let consider a principal ideal
$I(\overline{a})$, for $a \in Q$, of  ideal $\omega [Q]$ which is
not simple. Then, there  exists an element $b \in Q$ such that
$I(\overline{a}) = I(\overline{b}) + I(\overline{c})$, where $a =
bc$, and $I(\overline{b})$ is a proper ideal of algebra
$I(\overline{a})$}.\vspace*{0.1cm}

\textbf{Proof.} Let $J$ be a proper ideal of $I(\overline{a})$. By
Lemma 2.3 the normal subloops $B$ and $A$ of the loop $Q$
correspond to the ideals $J$, $I(\overline{a})$  and, by item 2)
of Lemma 4.6, $B \subset A$. Let $b \in B$ and let $a = bc$. Then
$b, c \in A$ and, by item 1) of Lemma 4.6, $e - b, e - c \in
\omega [A] \subseteq I(\overline{a})$. Hence $I(\overline{b})
\subseteq I(\overline{a})$, $ I(\overline{c}) \subseteq
I(\overline{a})$. By (35), $I(\overline{a}) \subseteq
I(\overline{b}) + I(\overline{c})$. Then $I(\overline{a}) =
I(\overline{b}) + I(\overline{c})$, as required.$\Box$

\vspace{3mm}

Let $F[Q] \in \mathcal{P}$ and let $I(\overline{a})$, where
$\overline{a} = e - a$, for $a \in Q$, be a principal ideal of
$F[Q]$. As $F[Q] \in \mathcal{P}$ then, by item 2) of Corollary
4.4, $I(\overline{a}) \in \mathcal{P}$. Let $A$ be the normal
subloop of the loop $Q$ induced, via Lemma 2.3, by the ideal
$I(\overline{a})$.  Then, by item 3) of Corollary 4.5,
$I(\overline{a}) = F[A]$ and any element of $I(\overline{a})$ has
the form
$$\sum_{i=1}^n\alpha_iu_i, \quad \sum_{i=1}^n\alpha_i \neq 0,\eqno{(36)}$$ where
$\alpha_i \in F$, $u_i \in Q$.

Further, for principal ideals $I(\overline{a_i})$,
$I(\overline{b_j}), \ldots$ we use the notations $F[A_i] =
I(\overline{a_i})$, $F[B_j] = I(\overline{b_j}), \ldots$ The
symbols $\overline{FY}$ $\overline{F[A]}$ will denote the
$F$-modules $FY$, $F[A]$.

\vspace{3mm}

\textbf{Lemma 5.4.} \textit{Let $F[Q] \in \mathcal{P}$ and let
$I(\overline{a})$, where $\overline{a} = e - a$, $a \in Q$, be a
principal ideal of $F[Q]$. Then, there exists an element $b \in Q$
such that $I(\overline{a}) = I(\overline{b}) + I(\overline{c})$,
where $a = bc$, $I(\overline{b})$ is a proper ideal of algebra
$I(\overline{a})$ and $\overline{F[A]} = \overline{F[B]} \oplus
M[K]$, where $M[K]$ denotes the $F$-submodule of $\overline{F[C]}$
generated by set $K = A \backslash B$}.\vspace*{0.1cm}

\textbf{Proof.} Let $Q = L/H$, where $L$ is a free Moufang loop.
We consider the homomorphisms $\varphi: LX \rightarrow LX/I =
F[X]$, $\psi: F[X] \rightarrow F[L]/\omega[H] = F[Q]$ (see
Proposition 2.4). From item 5) of Corollary 3.5, it follows that
any element in $\omega[H]$ has the form
$$\sum_{j=1}^m\beta_jh_j, \quad \sum_{j=1}^m\beta_j = 0, \eqno{(37)}$$
where $\beta_j \in F$, $h_j \in H$. By Lemma 5.3, $I(\overline{a})
= I(\overline{b}) + I(\overline{c})$ and $I(\overline{b})$ is a
proper ideal of $I(\overline{a})$.

If we denote $\psi^{-1}(A) = X_A$, $\psi^{-1}(B) = X_B$,
$\psi^{-1}(C) = X_C$, then \break $\psi^{-1}(I(\overline{a})) =
\psi^{-1}(F[A]) = F[X_A]$, $\psi^{-1}(I(\overline{b})) = F[X_B]$,
$\psi^{-1}(I(\overline{c})) = F[X_C]$. By (37), the homomorphism
$F[X] \rightarrow F[X]/\omega[H]$ maintains the sum of
coefficients, thus any element in $F[X_A]$, $F[X_B]$, $F[X_C]$ has
the form
$$\sum_{i=1}^k\gamma_ix_i, \quad \sum_{i=1}^k\gamma_i \neq 0, \eqno{(38)}$$
where $\gamma_i \in F, x_i \in X$. Then from (37),  (38) it
follows that $F[X_A] \bigcap \omega[H] = \{0\}$, $F[X_B] \bigcap
\omega[H] = \{0\}$, $F[X_C] \bigcap \omega[H] = \{0\}$.
Consequently, $F[A] = \psi(F[X_A]) = (F[X_A] +
\omega[H])/\omega[H] \cong F[X_A]/(F[X_A] \bigcap \omega[H]) =
F[X_A]/\{0\} \break = F[X_A]$, i.e. $F[X_A] \cong F[A]$.
Similarly, $F[X_B] \cong F[B]$, $F[X_C] \cong F[C]$.

According to  Lemma 2.1 $\varphi^{-1}(X_A) = X_A$,
$\varphi^{-1}(X_B) = X_B$, $\varphi^{-1}(X_C) = X_C$. Hence
$\varphi^{-1}(F[X_A]) = FX_A$, $\varphi^{-1}(F[X_B]) = FX_B$,
$\varphi^{-1}(F[X_C]) = FX_C$.

From the definition of an ideal $I$ of the loop algebra $FX$ it
follows that any  element of $I$ has the form
$\sum_{j=1}^n\beta_jx_j$ with $\sum_{j=1}\beta_j = 0$, where
$\beta_j \in F, x_j \in X$. Then from (38) it follows that $F[X_A]
\bigcap I = \{0\}$, $F[X_B] \bigcap I = \{0\}$, $F[X_C] \bigcap I
= \{0\}$ and  $FX[A] = \varphi(FX_A) = (FX_A + I)/I \cong \break
FX_A/(FX_A \bigcap I)  = F[X_A]/\{0\} = F[X_A]$. Hence $FX_A \cong
F[X_A]$. Before we have proved that $F[X_A] \cong F[A] =
I(\overline{a})$. Consequently, $FX_A \cong I(\overline{a})$.
Similarly, $FX_B \cong F[B] = I(\overline{b})$, $FX_C \cong F[C] =
I(\overline{c})$.

The inverse image of equality $I(\overline{a}) = I(\overline{b}) +
I(\overline{c})$ regarding  homomorphism $\varphi\psi$, is the
equality $FX_A = FX_B + FX_C$ of the loop algebra $FX$. The loop
algebra $FX$ is a free $F$-module with basis $\{x \vert x \in
X\}$. Then $\overline{FX_A} = \overline{FX_B} \oplus M(X_A
\backslash X_B)$, $M(X_A \backslash X_B) \subseteq
\overline{FX_C}$. Hence $\overline{F[A]} = \overline{F[B]} \bigcap
M[A \backslash B]$, $M[A \backslash B] \subseteq \overline{F[C]}$.
This completes the proof of Lemma 5.4.\vspace*{0.1cm}

\textbf{Proposition 5.5.} \textit{Let $F[Q] \in \mathcal{P}$ and
let $I(\overline{a})$, where $\overline{a} = e - a$, for $a \in
Q$, be a principal ideal of $F[Q]$. Then $I(\overline{a})$
decomposes into a direct sum of finite number of simple
nonassociative principal ideals   $I(\overline{a}) =
I(\overline{b}_1) \oplus \ldots \oplus
I(\overline{b}_n)$.}\vspace*{0.1cm}

\textbf{Proof.} Inductively we construct two series
$$I(\overline{b}_1) \supset I(\overline{b}_2) \supset \ldots
\supset I(\overline{b}_n) \supset \ldots ,$$
$$I(\overline{d}_1) \subseteq I(\overline{d}_2) \subseteq \ldots
\subseteq I(\overline{d}_n) \subseteq \ldots \eqno{(39)}$$ of
proper non-zero ideals of the algebra $I(\overline{a})$ such that
$I(\overline{a}) = I(\overline{b}_n) + I(\overline{d}_n)$ and a
series $$M[K_1] \subset M[K_2] \subset \ldots M[K_n] \subset
\ldots \eqno{(40})$$ of $F$-submodules of the $F$-module
$\overline{F[A]}$ such that $M[K_i] \subseteq M[D_i]$ and
$\overline{F[A]} = \overline{F[B_i]} \oplus M[K_i]$, $K_i = A
\backslash B_i$. The inductive process stops if an ideal
$I(\overline{b_n})$ is simple for some integer $n$.

Let the ideal $I(\overline{a})$ be non-simple. Then by Lemma 5.4
$I(\overline{a}) = I(\overline{b}_1) + I(\overline{c}_1)$, where
$I(\overline{b}_1)$ is a  proper ideal of $I(\overline{a})$ and
$\overline{F[A]} = \overline{F[B_1]} \oplus M[D_1]$, $D_1 = A
\backslash B_1$, $M[D_1] \subseteq \overline{F[C_1]}$. If at least
one of the ideals $I(\overline{b}_1)$, $I(\overline{c}_1)$ is
simple then the inductive process ends. Conversely, let us
consider that the ideal $I(\overline{b_1})$ is non-simple. By
Lemma 5.4, let $I(\overline{b}_1) = I(\overline{b}_2) +
I(\overline{c}_2)$, where $I(\overline{b}_2)$ is an ideal of
$F[A]$ and is a proper ideal of $I(\overline{b}_1)$. Again by
Lemma 5.4 $I(\overline{a}) = I(\overline{b}_2) +
I(\overline{d}_2)$ and $\overline{F[A]} = \overline{F[B_2]} \oplus
M[K_2]$, $K_2 = A \backslash B_2$, $M[K_1] \subseteq
\overline{F[D_2]}$.

Let us continue the inductive process. Let $I(\overline{a}) =
I(\overline{b}_n) + I(\overline{d}_n)$, $\overline{F[A]} =
\overline{F[B_n]} \oplus M[K_n]$, $K_n = A \backslash B_n$,
$M[K_n] \subseteq \overline{F[C_n]}$  and let the ideal
$I(\overline{b}_n)$ be non-simple. By Lemma 5.4 $I(\overline{b}_n)
= I(\overline{b}_{n+1}) + I(\overline{c}_{n+1})$ and
$I(\overline{a}) = I(\overline{b}_{n+1}) + I(\overline{d}_{n+1})$
and $\overline{F[A]} = \overline{F[B_{n+1}]} \oplus M[K_{n+1}]$,
$K_{n+1} = A \backslash B_{n+1}$.  From $I(\overline{b}_n) \supset
I(\overline{b}_{n+1})$ it follows that $I(\overline{d}_n) \subset
I(\overline{d}_{n+1})$ and $K_n \subset K_{n+1}$, $M[K_n] \subset
M[K_{n+1}]$. Consequently, the series (39), (40) with property
$\overline{F[A]} = \overline{F[B_n]} \oplus M[K_n]$, $n = 1, 2,
\ldots $ are defined.

The modules  $M[K_n]$ in the ascending series (40)  satisfy the
property $\overline{b}_1 \notin M[K_n]$. Then, by Zorn Lemma, this
series have a maximal proper ideal $J$ in $I(\overline{a})$ such
that $\overline{b}_1 \notin J$. Let $\overline{b} \in
I(\overline{a}) \setminus J$. As $J$ is a maximal ideal of
$I(\overline{a})$ then $I(\overline{a}) = I(\overline{b}) + J$.

Let the ideal $I(\overline{b})$ be non-simple. Then, by Lemma 5.4,
$I(\overline{b}) = I(\overline{b}_1) + I(\overline{b}_2)$, where
$I(\overline{b}_1)$ is a proper non-zero ideal of
$I(\overline{a})$, $I(\overline{a}) = I(\overline{b}_1) +
I(\overline{b}_2) + J$ and the ideal $I(\overline{b}_2) + J$
strictly contain the maximal ideal $J$. Contradiction. Hence the
ideal $I(\overline{b})$ is simple.

By Lemma 5.4 $I(\overline{a}) = I(\overline{b}) + I(\overline{d})$
for some proper ideal $I(\overline{d})$ of $I(\overline{a})$. The
ideal $I(\overline{b})$ is simple. Let the ideal $I(\overline{d})$
is non-simple. Then  $I(\overline{a}) = I(\overline{b}_1) +
I(\overline{d}_1)$, where $I(\overline{d}_1) = I(\overline{d})$
and $I(\overline{b}_1) = I(\overline{b})$ is a simple ideal.
Further, $I(\overline{a}) =  I(\overline{b}) + I(\overline{b}_2) +
\ldots + I(\overline{b}_k) + \ldots$, where $I(\overline{b}_i)$,
$i = 1, 2, \ldots$, is a simple ideal. Then, by Lemma 5.2,
$I(\overline{a}) =  I(\overline{b}) \oplus I(\overline{b}_2)
\oplus \ldots \oplus I(\overline{b}_k) \oplus \ldots$, where each
$I(\overline{b}_i)$, $i = 1, 2, \ldots$, is a simple ideal.

It is known  that any element of the direct sum is written
unequivocally as the sum of a finite number of  non-zero elements,
taken one from  some ideals $I(\overline{b}_{i_j})$. Let
$\overline{a} = \overline{b}_1 + \ldots + \overline{b}_n$, where
$\overline{b}_j \in I(\overline{b}_{i_j})$. Then by (34)
$I(\overline{a}) \subseteq I(\overline{b}_{i_1}) \oplus \ldots
\oplus I(\overline{b}_{i_n})$ and,  consequently, $I(\overline{a})
= I(\overline{b}_{i_1}) \oplus \ldots \oplus
I(\overline{b}_{i_n})$. As $F[Q] \in \mathcal{P}$ then by
\cite[Theorem 3, cap. 8]{ZSSS} $I(\overline{b_{i_j}}) \in
\mathcal{P}$ and by item 3) of Corollary 4.5 the ideals
$I(\overline{b_{i_j}})$ are nonassociative. This completes the
proof of Proposition 5.5.\vspace*{0.1cm}

\textbf{Lemma 5.6.} \textit{Let $Q \in \mathcal{T}$ be a
nonassociative semisimple Moufang loop, let $F[Q] \in \mathcal{P}$
be it corresponding alternative loop algebra and let $A, B$ be
normal subloops of the loop $Q$. Then $A \subset B$ (respect. $A =
B$) when and only when $\omega[A] \subset \omega[B]$ (respect.
$\omega[A] = \omega[B]$).}

\vspace{3mm}

\textbf{Proof.} Let $Q = L/H$, where $L$ is a free Moufang loop
and let $\varphi: LX \rightarrow LX/I = F[X]$, $\psi: F[X]
\rightarrow F[L]/\omega[H] = F[Q]$ be the homomorphisms considered
in proof of Lemma 5.4. Let $\psi^{-1}(A) = X_A$, $\psi^{-1}(B) =
X_B$, $X_A, X_B \subseteq X$. It is proved that $FX_A \cong F[X_A]
\cong F[A] = \omega[A]$, $FX_B \cong F[X_B] \cong F[B] =
\omega[B]$. From (36) -- (38), it follows easily  that the
restrictions of homomorphism $\psi$ on $X_A$ and on $X_B$ are
isomorphisms of loops $X_A, A$ and $X_B, B$, respectively.

From the mentioned isomorphisms, it follows that the inclusions
$F[A] \subset F[B]$ in the alternative loop algebra $F[Q]$ and
$FX_A \subset FX_B$ in the loop algebra $FX$ are equivalent. The
loop algebra $FX$ is a free $F$-module with basis $\{x \vert x \in
X\}$. Then the inclusion $FX_A \subset FX_B$ is equivalent to
inclusion $X_A \subset X_B$ of subloops in the loop $X$. Further,
from isomorphisms of loops $X_A, A$ and $X_B, B$ it follows that
the inclusion $X_A \subset X_B$ is equivalent to inclusion $A
\subset B$. Consequently, the inclusion $F[A] \subset F[B]$ in the
alternative loop algebra $F[Q]$ is equivalent to inclusion $A
\subset B$ of subloops in the loop $Q$. The facts that the
equalities $F[A] = F[B]$ and $A = B$ are equivalent are
analogously proved.

\vspace{3mm}

\textbf{Proposition 5.7.} \textit{Let $Q \in \mathcal{T}$ be
nonassociative semisimple Moufang loop and let $F[Q] \in
\mathcal{P}$ be its corresponding alternative loop algebra. Then
for any element $a \in Q$ the normal subloop $N(a)$ of $Q$
generated by element $a$ decompose into a direct product  of
finite number of nonassociative simple loops.}\vspace*{0.1cm}

\textbf{Proof.} By Proposition 5.5  $I(\overline{a}) =
I(\overline{a}_1) \oplus \ldots \oplus I(\overline{a}_k)$, where
each $I(\overline{a}_i)$ is a simple ideal of $F[Q]$ generated by
element $\overline{a}_i = e - a_i$, for $a_i \in Q$, $i = 1,
\ldots , k$. By Lemma 2.3, the ideal $I(\overline{a}_i)$ induces
in $Q$ the normal subloop $H_i = Q \bigcap (e +
I(\overline{a}_i))$. Let $N(a_i)$ denote the normal subloop of $Q$
generated by the element $a_i \in Q$. It is clear that $a_i \in
H_i$. Then $N(a_i) \subseteq H_i$. If $N(a_i) \subset H_i$
(strictly) then by Lemma 5.6 $\omega[N(a_i)] \subset \omega[H_i]$
(strictly). But $\omega[H_i] = I(\overline{a})$. Hence
$\omega[N(a_i)] \subset I(\overline{a})$ (strictly), i.e.
$\omega[N(a_i)]$  is a proper ideal of $I(\overline{a}_i)$. We get
a contradiction because $I(\overline{a}_i)$ is a simple ideal.
Consequently, $\omega[N(a_i)] = I(\overline{a}_i)$.

If $K$ is a proper normal subloop of $N(a_i)$ then, by Lemma 5.6,
$\omega[K]$ is a proper ideal of $\omega[N(a_i)] =
I(\overline{a}_i)$. Again we get a contradiction. Hence the normal
subloops $N(a_i)$, $i = 1, \ldots , k$,  are simple. We have

$$I(\overline{a}) = I(\overline{a}_1) \oplus \ldots \oplus
I(\overline{a}_k) \eqno{(41)}$$ or $\omega[N(a)] = \omega[N(a_1)]
+ \ldots + \omega[N(a_k)]$. Then, by item 2) of Lemma 4.6,
$\omega[N(a)] = \omega[N(a_1)]\cdot \ldots \cdot N(a_k)]$ and, by
Lemma 5.6, $N(a) = N(a_1)]\cdot \ldots \cdot N(a_k)$. The subloops
$N(a_i)$, $i = 1, \ldots , k$, are simple. Then, by Lemma 5.1,
$$N(a) = N(a_1) \otimes \ldots \otimes N(a_k). \eqno{(42)}$$ This completes the
proof of Proposition 5.7.$\Box$

\vspace{3mm}

\textbf{Corollary 5.8.} \textit{Let $Q \in \mathcal{T}$ be a
nonassociative semisimple Moufang loop and let $F[Q] \in
\mathcal{P}$ be it corresponding alternative loop algebra. Then:}

\textit{1) any nonassociative simple subloop of the loop $Q$ has
the form $H = I(\overline{a})$, where $I(\overline{a})$ is a
normal subloop of $F[Q]$, with $\overline{a} = e - a$ for some $a
\in Q$;}

\textit{2) any nonassociative simple subalgebra  of algebra $F[Q]$
has the form $F[H] = \omega[H]$, where $H = I(\overline{a})$,
with $\overline{a} = e - a$, for some $a \in Q$;}

\textit{3) a nonassociative subalgebra $F[H]$ of algebra $F[Q]$ is
simple when and only when the nonassociative normal subloop $H$ of
loop $Q$ is simple.}\vspace*{0.1cm}

The corollary follows from (41), (42) and Lemma
5.6.\vspace*{0.1cm}

\textbf{Lemma 5.9.} \textit{Any algebra $A$ of semisimple class
$\mathcal{P}$ of radical $\mathcal{R}$ decomposes into a direct
sum of nonassociative simple algebras.}\vspace*{0.1cm}

\textbf{Proof.} By item 3) of Corollary 4.5, any algebra $A$ of
semisimple class $\mathcal{P}$ has the form $A = F[Q]$ and $F[Q] =
\omega [Q]$. The ideal $\omega [Q]$ is generated as ideal by set
$\{e - g \vert g \in Q\}$. Let $g_1 \in Q$. As $(e - g)g_1 = (e -
gg_1) - (e - g_1)$, $g_1(e - g) = (e - g_1g) - (e - g_1)$ then
$\omega [Q]$ is generated as $F$-module by elements of form $e -
g$, where $g \in Q$. Denote $e - g = \overline{g}$ and let
$I(\overline{g})$ be the (principal) ideal generated by element
$\overline{g} \in F[Q]$. Then
$$F[Q] = \sum_{g \in Q}I(\overline{g}). \eqno{(43)}$$

As $F[Q] \in \mathcal{P}$ then by item 2) of Corollary 4.4
$I(\overline{g}) \in \mathcal{P}$. By Proposition 5.5,
$I(\overline{g}) = I(\overline{b}_1) \oplus \ldots \oplus
I(\overline{b}_k)$, where $I(\overline{b}_i)$ is a simple ideal of
$F[Q]$. Then, from (43), it follows that $F[Q] = \sum
I(\overline{b}_i)$  and, by Lemma 5.2, $F[Q] = \sum^{\oplus}
I(\overline{b}_i)$, where $I(\overline{b}_i)$ are  simple ideals
of algebra $F[Q]$. This completes the proof of Lemma
5.9.\vspace*{0.1cm}

\textbf{Corollary 5.10.} \textit{Let $char F = 0$ or $char F = 3$.
Then any nonassociative commutative Moufang loop $Q$ is
$\mathcal{S}$-radical and any it alternative loop algebra $F[Q]$
is $\mathcal{R}$-radical.}

\textbf{Proof.} According to Theorem 4.2 $F[Q]/\mathcal{R}(F[Q]) =
P(F[Q])$, where $P(F[Q]) \in \mathcal{P}$. We assume that $P(F[Q])
\neq Fe$. As the loop $Q$ is commutative then by (31) the algebra
$P(F[Q])$ also is commutative. From Lemma 5.9 it follows that
$P(F[Q])$ decomposes into a direct sum of nonassociative simple
algebras. But any commutative simple alternative algebra is a
field \cite[pag. 172]{ZSSS}. We get a contradiction. Hence
$P(F[Q]) = Fe$. Then, by using the definitions, $F[Q] \in
\mathcal{R}$ and $Q \in \mathcal{S}$, as required.$\Box$

\vspace{3mm} \textbf{Lemma 5.11.} \textit{Any nonassociative
semisimple Moufang loop $Q \in \mathcal{T}$ decomposes into a
direct product of nonassociative simple loops.}\vspace*{0.1cm}

The statement follows from (41) - (43) and Lemmas 5.1, 5.6.

\section{Main results}

Let us consider the following notions. In the beginning of the
paper, we have mentioned that, in the literature, an algebra is
called antisimple, if none of its two-sided ideals allows
homomorphism on a simple algebra.

An alternative loop algebra $F[Q] \in \mathcal{M}$ will be called
\textit{antisimple with respect to nonassociativity} if for any
its ideal $\omega[H]$ the algebra $(\omega[H])^{\sharp} = F[H]$
does not allow homomorphism on a simple nonassociative algebra.

Analogously, a loop $Q \in \mathcal{L}$ will be  called
\textit{antisimple with respect to  nonassociativity}  if none of
its normal subloops allows homomorphism on a simple nonassociative
loop.

Let $\mathcal{A}$ be a class of algebras, let $\mathcal{B}$ be its
radical class and let $\mathcal{C}$ be it semisimple class. By the
definition of radical $\mathcal{B}$ any homomorphic image of
$\mathcal{B}$-algebra is a $\mathcal{B}$-algebra. In \cite[
Proposition 1, pag. 184]{ZSSS} it is proved that the radical class
$\mathcal{B}$ of $\mathcal{A}$ is the totality of algebras from
 $\mathcal{A}$, not reflected homomorphically on the algebras of class
$\mathcal{C}$. Then from Lemmas 5.9, 5.11  it follows the next
result.

\vspace{3mm}

\textbf{Lemma 6.1.} \textit{The class of all antisimple with
respect to nonassociativity alternative loop algebras $F[Q]$
coincides with the  radical class $\mathcal{R}_A$ of all
alternative loop algebras of type $F[Q] = (\omega[Q])^{\sharp}$.
The class of all antisimple with respect to  nonassociativity
Moufang  loops coincides with the radical class $\mathcal{S}$ of
Moufang loops.}\vspace*{0.1cm}

\textbf{Proposition 6.2} \textit{An alternative loop algebra
$F[Q]$ is antisimple with respect to  nonassociativity when and
only when $F[Q]$ does not have subalgebras that are nonassociative
simple algebras. A Moufang loop $Q$ is antisimple with respect to
nonassociativity when and only when  $Q$ does not have subloops
that are nonassociative  simple loops.}\vspace*{0.1cm}

\textbf{Proof.} Let $F[Q] \in \mathcal{M}$ and, according to
Theorem 4.2, let $F[Q]/\mathcal{R}(F[Q]) \break = P(F[Q])$, where
$P(F[Q]) \in \mathcal{P}$. By Lemma 6.1, the algebra $F[Q]$ is
antisimple with respect to nonassociativity when and only when
$P(F[Q]) = Fe$.

We assume that $P(F[Q]) \neq Fe$. By Lemma 2.3, the ideal
$\mathcal{R}(F[Q])$ of algebra $F[Q]$ induces the normal subloop
$R = Q \bigcap (e +  \mathcal{R}(F[Q]))$ of loop $Q$. As $P(F[Q])
\neq Fe$ then $\mathcal{R}(F[Q])$ is a proper ideal of $F[Q]$ and,
according to the one-to-one mapping among all ideals of  $F[Q]$
and all normal subloops of $Q$ the normal subloop $R$ is proper
and $F[Q/R] \cong F[Q]/\mathcal{R}(F[Q]) = P(F[Q])$. By Lemma 5.9,
$P(F[Q]) \neq Fe$ decomposes into a direct sum of nonassociative
simple algebras. According to Corollary 5.8, let $F[\overline{H}]
= \omega[\overline{H}]$, where $\overline H$ is the normal subloop
of loop $Q/R$ generated by one element $hR = h +
\mathcal{R}(F[Q])$, be one of such nonassociative simple algebras.
We denote by $H$ the normal subloop of the loop $Q$ generated by
the element $h \in Q$.

Clearly,  the inverse image of subalgebra $F[\overline{H}]$ under
the natural homomorphism $F[Q] \rightarrow F[Q]/\mathcal{R}(F[Q])
= P(F[Q])$ is $F[H] + \mathcal{R}(F[Q]$.  For $a \in
\mathcal{R}(F[Q])$, we have $a = \sum_{q \in Q}\alpha_qq$ with
$\sum_{q \in Q}\alpha_q = 0$, by item 5) of Corollary  3.5, and
for $b \in F[\overline{H}]$, we have $b = \sum_{\overline{g} \in
\overline{H}}\beta_{\overline{g}}\overline{g}$ with
$\sum_{\overline{g} \in \overline{H}}\beta_{\overline{g}} \neq 0$,
by item 2) of Corollary 4.5. The extension
$F[Q]/\mathcal{R}(F[Q])$ does not change the sum of coefficients.
Hence if $c \in F[H]$ then $c = \sum_{g \in H}\gamma_gg$ with
$\sum_{g \in H}\gamma_g \neq 0$. Then $F[H] \bigcap
\mathcal{R}(F[Q]) = (0)$ and, by homomorphism theorems, it follows
$F[\overline{H}] \cong (F[H] +
\mathcal{R}(F[Q]))/\mathcal{R}(F[Q]) \cong F[H]/(F[H] \bigcap
\mathcal{R}(F[Q]) = F[H]$. Hence, the subalgebra $F[H]$ of the
algebra $F[Q]$ is a nonassociative simple algebra. Then, by item
3) Corollary 5.8, $H$ is a nonassociative simple loop.
Consequently, if $P(F[Q]) \neq (0)$ then:

(i) the algebra $F[Q]$ contains a nonassociative simple algebra
$F[H]$;

(ii) the loop $Q$ contains a nonassociative simple loop $H$.

If the algebra $F[Q]$ does not contain a nonassociative simple
loop then, from item 2) of Corollary 4.5, it follows that $\sum_{q
\in Q}\alpha_q = 0$ for any element $a = \sum_{q \in Q}\alpha_qq
\in \omega[Q]$. In such a case, $F[Q] = (\omega[Q])^{\sharp}$ and,
from Lemma 4.1, it follows that $P(F[Q]) = Fe$. This completes the
proof of the first statement.

Now, let $Q \in \mathcal{L}$ be a Moufang loop and, according to
Theorem 4.7, let $Q/\mathcal{S}(Q) = T(Q)$, where $T(Q) \in
\mathcal{T}$. By Lemma 6.1, the loop $Q$ is antisimple with
respect to nonassociativity when and only when $T(Q) = \{e\}$. We
prove that the equality $T(Q) = \{e\}$ is equivalent to the
property that the loop $Q$ does not contain a subloop isomorphic
to simple nonassociative loop.

Indeed, we assume that $T(Q) \neq \{1\}$. Then, from the relation
$Q/\mathcal{S}(Q) = T(Q)$, it follows that $Q \neq \mathcal{S}(Q)$
and, from the definition of the class $\mathcal{S}$, it follows
that $\mathcal{R}(F[Q]) \neq F[Q]$ and $P(F[Q]) \neq Fe$, in
accordance with the relation \break $F[Q]/\mathcal{R}(F[Q])  =
P(F[Q])$. In such a case, the loop $Q$ contains a nonassociative
simple loop $H$, by statement (ii).

Now let the loop $Q$  not contain  any nonassociative simple loop
as sub\-loop. Then, from item 3) of Corollary 5.8, it follows that
the alternative loop algebra $F[Q]$ does not contains
nonassociative simple algebra as subalgebra. Thus, by the first
case, $P(F[Q]) = Fe$ or $\mathcal{R}(F[Q]) = F[Q]$,
$\mathcal{S}(Q) = Q$, $T(Q) = \{e\}$. This completes the proof of
Proposition 6.2.\vspace*{0.1cm}

Let us consider the  analogue for alternative loop algebras of the
Wedderburn Theorem for finite dimensional associative algebras.

By Kleinfeld Theorem \cite{ZSSS} any nonassociative simple
alternative algebra is a Cayly-Dickson algebra over its centre.
Then, from
 Theorem 4.2 and Lemmas 5.9, 6.1, it follows the next result.

 \vspace{3mm}

\textbf{Proposition 6.3.} \textit{Let $F[Q]$ be an alternative
loop algebra from class $\mathcal{M}$ and let $\mathcal{R}(F[Q])$
be its radical. Then  algebra $(\mathcal{R}(F[Q]))^{\sharp} =
F[G]$, $G \subseteq Q$, $F[G] \in \mathcal{R}_A$, is a
nonassociative antisimple with respect to nonassociativity or,
equivalently, it does not contain subalgebras that are
nonassociative simple algebras and the quotient-algebra
$F[Q]/\mathcal{R}(F[Q])$ is a direct sum of Cayley-Dickson
algebras over their centre.}\vspace*{0.1cm}

As it was above mentioned, the nonassociative antisimple with
respect to nonassociativity algebras are considered in
Propositions 3.3, 3.4, 6.2 and Corollary 3.5.

Let now $Q \in \mathcal{L}$ be a nonassociative Moufang loop.
According to Theorem 4.7, $Q/\mathcal{S}(Q) = T(Q)$, where $T(Q)
\in \mathcal{T}$, $\mathcal{S}$ is the radical class,
$\mathcal{T}$ is the semisimple class for class loop
$\mathcal{L}$. Further, as a rule in the  theory of algebraic
systems, in order to study the  loops of class $\mathcal{L}$ we
will consider the loops of classes $\mathcal{S}$ and $\mathcal{T}$
separately.

To describe  class $\mathcal{T}$, we remind the description of
nonassociative simple Moufang loops from \cite{San11}. Let $M(F)$
denote the \textit{matrix Paige loop} constructed, over the field
$F$, as in \cite{Paige}. That is, $M(F)$ consists of vector
matrices

\noindent $M^{\ast}(F) = \left(
\begin{array}{ll} \alpha_1 & \alpha_{12}\\
\alpha_{21} & \alpha_2  \end{array} \right)$, where $\alpha_1,
\alpha_2 \in F$, $\alpha_{12}, \alpha_{21} \in F^3$, \break
$\text{det}$ $M^{\ast} = \alpha_1\alpha_2 - (\alpha_{12},
\alpha_{21}) = 1$, and where $M^{\ast}$ is identified with
$-M^{\ast}$.

\noindent The multiplication in $M(F)$ coincides with the Zorn
matrix multiplication $$ \left(
\begin{array}{ll} \alpha_1 & \alpha_{12}\\
\alpha_{21} & \alpha_2  \end{array} \right) \left(
\begin{array}{ll} \beta_1 & \beta_{12} \\ \beta_{21} & \beta_2
\end{array} \right) = $$
$$ \left( \begin{array}{ll} \alpha_1 \beta_1 + (\alpha_{12},
\beta_{21}) & \alpha_1 \beta_{12} + \beta_2 \alpha_{12} -
\alpha_{21} \times \beta_{21} \\ \beta_1 \alpha_{21} + \alpha_2
\beta_{21} + \alpha_{12} \times \beta_{12} & \alpha_2 \beta_2 +
(\alpha_{21}, \beta_{12}) \end{array} \right), $$ where, for
vectors $\gamma = (\gamma_1, \gamma_2, \gamma_3)$, $\delta =
(\delta_1, \delta_2, \delta_3) \in A^3$,\   $(\gamma, \delta) =
\gamma_1\delta_1 + \gamma_2\delta_2 + \gamma_3\delta_3$ denotes
their scalar product and $\gamma \times \delta = (\gamma_2\delta_3
- \gamma_3\delta_2, \gamma_3\delta_1 - \gamma_1\delta_3,
\gamma_1\delta_2 - \gamma_2\delta_1)$ denotes the cross vector
product.

Let $\Delta$ be a prime field and let $P$ be its algebraic
closure. In \cite{San11}, it was proved that only and only the
Paige loops $M(F)$, where $F$ is a Galois extension over $\Delta$
in $P$ are, up to an isomorphism, nonassociative simple Moufang
loops. \cite{San11} also  describes the finite nonassociative
simple Moufang loops and  the set of generators, the group of
automorphisms of nonassociative simple Moufang loops.

From Theorem 4.7, Lemmas 5.10, 6.1 and Proposition 6.2 it follows
for Moufang loops an analogue of the Wedderburn Theorem, for
finite dimensional associative algebras.

\vspace{3mm}

\textbf{Theorem 6.4.} \textit{Let $\Delta$ be a prime field, $P$
be its algebraic closure, and $F$ be a Galois extension over
$\Delta$ in $P$. Then the radical $\mathcal{S}(Q)$ of a Moufang
loop $Q$ is nonassociative antisimple with respect to
nonassociativity or, equivalently, it does not contain any
subloops that are nonassociative simple loops and quotient-loop
$Q/\mathcal{S}(Q)$ is isomorphic to the direct product of matrix
Paige loops $M(F)$.}\vspace*{0.1cm}

Let $G$ be a finite Moufang loop. Obviously, from the finiteness
of $G$ it follows that for any subloop $H$ of $G$ there exists a
normal subloop $K$ of $H$ such that the composition factor $H/K$
is a simple loop. According to \cite{GrZ}, a finite Moufang loop
$G$ it said to be a loop of \textit{group type} if all composition
factors of $G$ are groups. It is clear that the homomorphic image
of loop of group type is a loop of group type and the product of
two normal subloops of group type is again a loop of group type.
Hence every finite Moufang loop has a unique maximal normal
subloop of group type (\cite[Proposition 1]{GrZ}). We denote this
maximal normal subloop of group type by $Gr(G)$. It is obvious
that $Gr(G/Gr(G)) = e$. Hence $Gr(G)$ is a radical of $G$. By
\cite{GrZ} $Gr(G)$ is called the \textit{group-type radical} of
$G$. Obviously, for any finite Moufang loop $G$, the definition of
normal subloop $Gr(G)$ is equivalent with the condition: $Gr(G)$
does not contain any subloops that are nonassociative simple
loops, by Theorem 6.4. Hence, for finite loops, the group-type
radical $Gr$ coincides with radical $\mathcal{S}$.

In the proof of the main result from \cite{GrZ} about the
existence of quasi-$p$-Sylow subloops in every finite Moufang
loop, the following structural Theorem B is used essentially:
\emph{every finite Moufang loop $M$ contains uniquely determined
normal subloops $Gr(M)$ and $M_0$ such  that $Gr(M) \leq M_0$,
$M/M_0$ is an elementary abelian $2$-group, $M_0/Gr(M)$ is the
direct product of simple Paige loops $M(q)$ (where $q$ may vary),
the composition factors of $Gr(M)$ are groups and} $Gr(M/Gr(M)) =
1$.

The proof of Theorem B is based on the correspondence between
Moufang loops and groups with triality \cite{Doro}. The proof of
Theorem B is quite cumbersome and uses deep results from finite
groups. Moreover, the Theorem B in such a version does not hold
true. For example, the case when $M$ is a simple loop leads to a
contradiction with the  condition that $M/M_0$ is an elementary
abelian $2$-group. In reality, $M/M_0$ is the unitary group, i.e.
$M = M_0$. In such a case Theorem B is a particular case  of
Theorem 6.4. Hence if the corresponding results from this paper
are used in the proofs of the main results from \cite{GrZ}, then
these  proofs become as it is shown below.

From Theorem 6.4, it follows  that the loops from semisimple class
$\mathcal{T}$ are well described. However, unlike the class
$\mathcal{T}$, much less is known about the qualities and
construction of loops from the radical class $\mathcal{S}$. A new
approach is suggested for the study of the loops in the class
$\mathcal{S}$ (though some authors made some attempts earlier):

a) by using the one-to-one mapping between loops $Q \in
\mathcal{S}$ and alternative loop algebras $F[Q] \in \mathcal{R}$,
below indicated in Theorem 6.5;

b) by using the developed  theory of alternative algebras, in
particular, of the algebras  with  externally adjoined unit, of
Zhevlakov radicals, of circle loops and others.

According to Lemma 6.1 and Proposition 6.2 the following
statements are equivalent for a nonassociative Moufang loop $Q \in
\mathcal{L}$:

r1) $Q \in \mathcal{S}$;

r2) $Q$ is an loop antisimple with respect to nonassociativity;

r3) the loop $Q$ does not have any subloops that are simple loops.

\vspace{3mm} Then the opposite statements

nr1) $G \notin \mathcal{S}$, i.e. $G \in \mathcal{L} \backslash
\mathcal{S}$;

nr2) $G$ is not antisimple with respect to nonassociativity loop;

nr3) the loop $G$ has subloops that are simple loops

\noindent hold for any nonassociative Moufang loop $G \in
\mathcal{L} \backslash \mathcal{S}$.

From the definition of class of alternative loop algebras
$\mathcal{R}_A$, definition of class of loops $\mathcal{S}$ and
Theorem 3.2, it follows that if a nonassociative Moufang loop $Q$
satisfies the condition r1) then the loop $Q$ can be embedded into
the loop of invertible elements $U(F[Q])$ of the alternative loop
algebra $F[Q]$. On the other hand, \cite{San08}  proves that if a
nonassociative  Moufang loop $G$ satisfies the condition nr2) then
the loop $Q$ is not imbedded into the loop of invertible elements
$\mathcal{U}(A)$ for a suitable unital alternative $F$-algebra
$A$, where $F$ is an associative commutative ring with unit. As
$\mathcal{S} \bigcap (\mathcal{L} \backslash \mathcal{S}) =
\varnothing$ then the main result of this paper follows from the
above statements.\vspace*{0.1cm}

\textbf{Theorem 6.5.} \textit{Any nonassociative Moufang loop $Q$
that satisfies one of the equivalent conditions r1) - r3) can be
embedded into a loop of invertible elements $U(F[Q])$ of
alternative loop algebra $F[Q]$. The remaining loops of class of
all nonassociative Moufang loops $\mathcal{L}$, i.e.the loops $G
\in \mathcal{L}$ that satisfy one of the equivalent conditions
nr1) - nr3) cannot be embedded into a loop of invertible elements
of any unital alternative algebras.}\vspace*{0.1cm}

From Corollary 5.10 and Theorem 6.5 the following Corollary
follows.\vspace*{0.1cm}

\textbf{Corollary 6.6.} \textit{Any commutative Moufang loop $Q$
can be embedded into a loop of invertible elements $U(F[Q])$ of
the commutative alternative loop algebra $F[Q]$.}\vspace*{0.1cm}

\textbf{Corollary 6.7.} \textit{Any finite Moufang $p$-loop $Q$
can be embedded into a loop of invertible elements $U(F[Q])$ of
alternative loop algebra $F[Q]$.}\vspace*{0.1cm}

\textbf{Proof.} According to \cite{San11}, with up to an
isomorphism, only Paige loops $M(q)$ over a finite field $F_q$ are
finite simple Moufang loops. By \cite{Paige} the order of $M(Q)$
is $(1/d)q^3(q^4 - 1)$, where $d = gcd(2, q - 1)$, and is the
product of two coprime numbers $q^3(q^2 - 1)$ and $(1/d)q^3(q^2 +
1)$. $Q$ is a finite $p$-loop if the order of elements from $Q$ is
a power of $p$ (for Moufang loops, this is equivalent to the
condition that the order of $Q$ be a power of $p$). From here it
follows that the finite $p$-loop satisfies the condition r3). By
Theorem 6.5 the Corollary 6.7 is proved.\vspace*{0.1cm}

From Theorem 6.5, Corollary 6.6 and (10) the following Corollary
follows.

\vspace{3mm}

\textbf{Corollary 6.8.} \textit{Any nonassociative Moufang loop
$Q$ that satisfy one of the equivalent  conditions r1) - r3) can
be embedded into a circle loop $(U^{\star}(F[Q]), \circ)$ of
alternative loop algebra $F[Q]$. The remaining loops of class of
all nonassociative Moufang loops $\mathcal{L}$, i.e.the loops $G
\in \mathcal{L}$ that satisfy one of the equivalent conditions
nr1) - nr3) cannot be embedded into  circle loops of
 any unital alternative algebras.}

 \vspace{3mm}

\textbf{Corollary 6.9.} \textit{Any commutative Moufang loop $Q$
can be embedded into a circle loop   $(U^{\star}(F[Q]), \circ)$ of
the alternative loop algebra $F[Q]$.}

\vspace{3mm} In \cite{Good87}, the circle loops are examined under
condition that the underlying sets defining the alternative
algebra and the loop are identical. Under this supposition it is
proved that for any prime $p$, there are no Moufang circle loops
of order dividing $p^4$ which are not associative, though non
associative Moufang loops of order $p^4$ exist \cite{Chein}. It is
also proved that no commutative Moufang which is not associative
is the circle loop of an alternative nil ring of index 2. But from
Corollaries 6.6, 6.7 the following Corollary
follows.\vspace*{0.1cm}

\textbf{Corollary 6.10.} \textit{Any commutative Moufang loop $Q$
can be embedded into a circle loop $(U(F[Q], \circ, 0)$  of the
alternative loop algebra $F[Q]$.}

\vspace{3mm}
 \textbf{Corollary 6.11.}
\textit{Any finite Moufang $p$-loop $Q$ can be embedded into a
circle loop of invertible elements $(U(F[Q], \circ, 0)$ of the
alternative loop algebra $F[Q]$.}

\vspace{3mm}

Now let us present some  examples that were proved on the  basis
on the correspondence between commutative Moufang loops and loops
of invertible elements of commutative alternative algebras
(Corollary 6.6).\vspace*{0.1cm}

$1^0$. \textbf{The Bruck's Theorem.} This theorem is one of the
profound results in the theory of commutative Moufang loops:\emph{
a commutative Moufang loop with $n$ ($n \geq 2$) generators is
centrally nilpotent of class at most} $n - 1$ \cite{Bruck}, Chap.
VIII]. The proof of this assertion is very cumbersome; it is based
on a complicated inductive process and uses several hundred
nonassociative identities. In \cite[Chap. 1]{Manin}  Manin used
group methods to prove a weaker assertion, namely that any finite
commutative Moufang loop of period 3 is centrally nilpotent.
Although less calculative, his proof is by no ???means simple; it
uses deep facts from finite group theory. The supremum of the
central nilpotence class  of a commutative Moufang loop with $n$
generators is equal to $n - 1$ \cite{Malb}.

In \cite{San99}, the relationship between commutative Moufang
loops and alternative commutative algebras, i.e. the Corollary
6.6, (in \cite{San99} the proof is not very convincing) is used to
prove (rather simply) that any finitely generated commutative
Moufang loop is centrally nilpotent. In the proof, we use the fact
that any alternative commutative nil-algebra of index 3 is locally
nilpotent, only.

In \cite{Pchel}, a Moufang loop $E$ is called \emph{special}, if
it can be embedded in the loop  $U(A)$ of invertible elements of
an alternative algebra $A$ with unit. The  Bruck Theorem is proved
in a quite transparent manner (and an accurate appraisal is made)
for special commutative Moufang loops. In the proof  the assertion
that the commutator ideal of the multiplication algebra of a free
commutative alternative algebra with $n$ free generators is
nilpotent of index $n - 1$ is transferred on such loops.
Consequently, according to Corollary 6.5, in \cite{Pchel} the
Bruck Theorem (and the accurate appraisal) is proved for any
commutative Moufang loops.\vspace*{0.1cm}

$2^0$. \textbf{Infinite independent system of identities.} In
\cite{Dnies}  Slin'ko has formulated the question: if any variety
of solvable alternative algebras would be finitely based. Umirbaev
has got an affirmative answer to this question for alternative
algebras over a field of characteristic $\neq 2, 3$ (see
\cite{Umir}), while  Medvedev \cite{Med} gave a negative answer
for characteristic 2. The topic of work \cite{San04} is the
transfer of infinite independent systems of a commutative Moufang
loop, constructed in \cite{San87} on solvable alternative
commutative algebra over a field of characteristic 3 (another
example was constructed by  Badeev, see \cite{Bad}), provided that
holds the Corollary 6.6. Consequently, the last result together
with the former results, completes the statement of Slin'ko
problem for solvable alternative algebras.\vspace*{0.1cm}

$3^0$. \textbf{The order of free commutative Moufang loops of
exponent 3.} Let $L_n$ be the free commutative Moufang loop on $n$
generators of exponent $3$ with unit $e$ and let $\mid L_n\mid$ $=
3^{\delta(n)}$. The Manin problem asks to calculate $\delta(n)$
\cite{Manin}. One of the main results of paper \cite{GrShest09} is
that $\delta(3) = 4$, $\delta(4) = 12$, $\delta(5) = 49$,
$\delta(6) = 220$, $\delta(7) = 1014$ or $1035$ and $\delta(7) =
1014$ if and only if $L_n$ can be embedded into a loop of
invertible elements of a unital alternative commutative algebra.
 It is also  proved that the free loop $L_n$ on $n < 7$ generators
is embedded into a loop of invertible elements $U(A)$ for a unital
alternative commutative algebra $A$. Moreover, $L_7$ may be
embedded in $U(A)$ if and only if the following identity is true
for commutative Moufang loops
$$((((a,x,y),z,b),t,c),b,c)((((a,x,z),y,b),t,c),b,c)$$
$$((((a,x,t),y,b),z,c),b,c)^{-1}((((a,x,b),y,z),t,c),b,c)
\eqno{(44)}$$
$$((((a,x,c),y,z),t,b),b,c)((((a,x,b),y,c),z,t),b,c) = e.$$
Here $(x,y,z) = (xy)z \cdot (x(yz))^{-1}$. According to Corollary
6.6, from here it follows that $\delta(7) = 1014$ and the identity
(44) is true for any commutative Moufang loop.

\section{Finite Moufang $p$-loops}

Let $Q$ be a loop with unit $e$. The set $\{z \in Q\ \vert\ zx =
xz,\ zx\cdot y = z \cdot xy,\ xz \cdot y = x \cdot zy,\ xy \cdot z
= x \cdot yz \quad \forall x, y \in Q \}$ is a  subloop $Z(Q)$ of
$Q$, the \textit{centre}. $Z(Q)$ is an abelian group, and every
subgroup of $Z(Q)$ is a normal subloop of $Q$. If $Z_1(Q) = Z(Q)$,
then the normal subloops $Z_{i + 1}(Q): Z_{i + 1}(Q)/Z_i(Q) =
Z(Q/Z_i(Q))$ are inductively determined. A loop $Q$ is called
\textit{centrally nilpotent of class $n$}, if its \textit{upper
central series} have the form $$\{e\} \subset Z_1(Q) \subset
\ldots \subset Z_{n - 1} \subset Z_n(Q) = Q.$$ If $N$ is a normal
subloop of $Q$, there is a unique smallest normal subloop $M$ of
$Q$ such that $N/M$ is part of the centre of $Q/M$, and we write
$M = [N,Q]$. The \textit{ lower central series} of $Q$ is defined
by $Q_1 = Q, Q_{i + 1} = [Q_i,Q]$ ($i \geq 1$). The loop $Q$ is
centrally nilpotent of class $n$ if and only if its lower central
series have the form $Q \supset Q_1 \supset \ldots \supset Q_{n -
1} \supset Q_n = \{e\}$ \cite{Bruck}.

The \textit{associator} $(x,y,z)$ and \textit{commutator} $(x,y)$
of elements $x, y, z \in Q$ are defined by the equalities $xy\cdot
z = (x\cdot yz)(x,y,z)$ and $xy = (yx)(x,y)$ for an arbitrary loop
$Q$. The \textit{commutator-associator of weight $n$} is defined
inductively:

1) any associator $(x,y,z)$ and any commutator $(x,y)$, where $x,
y, z \in Q$, are commutator-associator  of the weight $1$;

2) if $a$ is a commutator-associator  of weight $i$, then
$(a,x,y)$ or $(a,x)$, where $x, y \in Q$, is a
commutator-associator  of the weight $i+1$.\vspace*{0.1cm}

\textbf{Lemma 7.1} \cite{CovSan}. \textit{The subloops $Q_{i}$ ($i
= 1, 2, \ldots$) of the lower central series of a Moufang loop $Q$
are generated by all commutator-associators  of weight $i$ of
$Q$.}

\vspace{3mm}

If $A$ is an $F$-algebra, then its $n$ degree $A^n$ is an
$F$-module with a basis, consisting of products from any of its
$n$ elements with any bracket distribution. Algebra $A$ is called
\textit{nilpotent} if $A^n = (0)$ for a certain $n$.

\vspace{3mm}

\textbf{Lemma 7.2.} \textit{Let $Q$ be a finite Moufang $p$-loop
and $F$ be a field  of characteristic $p$. Then, the ideal
$\omega[Q]$ of the alternative loop algebra $F[Q]$ is
nilpotent.}\vspace*{0.1cm}

\textbf{Proof.}  In accordance with Corollary 6.7 we assume that
$Q \subseteq F[Q]$. By Theorem 1.2 from \cite[pag. 92]{Bruck} in a
finite Moufang loop $Q$ the order of any of its element divides
the order of $Q$. Hence  $g^k = e$, where $k = p^n$, for $g \in
Q$. We have $(e - g)^k = e - C_k^1g + \ldots + (-1)^iC_k^ig^i +
... + (-1)^kg^k$. All binomial coefficients $C_k^i$ can be divided
by $p$, therefore $(e - g)^k = e + (-1)^kg^k$. If $p = 2$, then
$(e - g )^k = e + g^k = e + e = 2 = 0$, because $F$ is a field of
characteristic 2. But if $p > 2$, then $(e - g)^k = e - g^k = e -
e = 0$. Then we can apply the following statement to algebra
$\omega Q$: any alternative $F$-algebra, generated as an
$F$-module by a finite set of nilpotent elements, is nilpotent
\cite[pages 144, 408]{ZSSS}. Consequently the augmentation ideal
$\omega[Q]$ is nilpotent, as required.\vspace*{0.1cm}

Let now $A$ be an alternative $F$-algebra with  unit $e$ and $B$
be a subalgebra from $A$, satisfying the law
$$x^m = 0. \eqno{(45)}$$
Then $e - B = \{e - b \vert b \in B\}$ will be a loop and $(e -
b)^{-1} = e + b + \ldots + b^{m-1}$.  We remind that inscription
$(a,b,c) = ab\cdot c - a\cdot bc, (a,b) = ab - ba$ means the
associator and commutator in algebra, but $[a,b,c] = (a\cdot
bc)^{-1}\cdot (ab\cdot c), [a,b] = a^{-1}b^{-1}\cdot ab$ are
\emph{associator} and \emph{commutator} in $IP$-loops.

\vspace{3mm}

\textbf{Lemma 7.4.} \textit{Let $A$ be an alternative algebra with
unit $e$ and $B$ its subalgebra, satisfying the law (45). Then for
$u, v, w \in B$ $[e - u , e - v, e - w] = e - ((e + w + \ldots +
w^{m-1})(e + v + \ldots + v^{m-1})\cdot(e + u + \ldots +
u^{m-1}))(u,v,w), [u,v] = (e + u + \ldots
 + u^{m-1})(e + v + \ldots + v^{m-1})(u,v)$.}
\vspace*{0.1cm}

\textbf{Proof.}  We denote $e - u = a,\ e - v = b,\ e - w = c$.
Then we have $[e - u, e - v, e - w] = (a\cdot bc)^{-1}(ab\cdot c)
= (a\cdot bc)^{-1}(ab\cdot c) - (a\cdot bc)^{-1}(a\cdot bc) + e =
e + (a\cdot bc)^{-1}(a,b,c) = e + (((e - w)^{-1}\cdot (e -
v)^{-1})(e - u)^{-1})(e - u, e - v, e - w) = e - ((e - w)^{-1}(e -
v)^{-1}\cdot (e - u)^{-1})(u,v,w) = e - ((e + w + \ldots +
w^{m-1})(e + v + \ldots + v^{m-1})\cdot (e + u + \ldots +
u^{m-1}))(u,v,w)$.  The  second  equality  is analogously  proved.

\vspace{3mm}

\textbf{Lemma 7.5.} \textit{Let $Q$ be a Moufang loop and let the
ideal $\omega[Q]$ of the alternative loop algebra $F[Q]$ be
nilpotent. Then the loop $Q$ is centrally nilpotent.}
\vspace*{0.1cm}

\textbf{Proof.}  It follows, from  the definition of  ideal
$\omega[Q]$, that $Q = Q_0  \subseteq e - \omega[Q] = e -
\omega[Q]^{2\cdot 0 + 1}$. We suppose that $Q_{i-1} \subseteq e -
(\omega[Q])^{2\cdot i + 1}$. Then it follows from  Lemma 7.4 that
$Q_i \subseteq e - (\omega[Q])^{2\cdot (i+1) + 1}$. The algebra
$\omega[Q]$ is nilpotent; we suppose that $(\omega[Q])^{2\cdot
k+3} = (0)$. Then $Q_k = e$. Hence the loop $Q$ is centrally
nilpotent, as required. $\Box$

\vspace{3mm} The following Proposition  follows from Lemmas 7.2
and 7.5.

\vspace{3mm}

\textbf{Proposition 7.6.} \textit{Any finite Moufang $p$-loop is
centrally nilpotent.}\vspace*{0.1cm}

\textbf{Theorem 7.7.} \textit{Let a Moufang loop $Q$  belong to
radical class $\mathcal{S}$. Then the following statements are
equivalent:}

\textit{1) the augmentation ideal $\omega[Q]$ of alternative loop
algebra $F[Q]$ is nilpotent;}

\textit{2) $Q$  is a finite $p$-loop and the field $F$ has a
characteristic $p$;}

\textit{3) the algebra $\omega Q$ is artinian, i.e. satisfies the
minimum condition for its left ideals.}\vspace*{0.1cm}

\textbf{Proof.}   Let the algebra $\omega[Q]$ be nilpotent, for
example, of index $r$ and let $0 \neq x \in (\omega[Q])^{r - 1}$.
By (31), the element $x$ will be written in form $x = \alpha_1g_1
+ \dots + \alpha_kg_k$, where $\alpha_i \in F$, $g_i \in Q$. We
suppose that $g_i \neq g_j$ if $i \neq j$. As $Q \in \mathcal{S}$
then from definition of the class $\mathcal{S}$ it follows that
$F[Q] \in \mathcal{R}$ and by item 4) of Proposition 3.3, we have
$x(e - u) = 0$ for any $u \in Q$. Hence $\alpha_1g_1 + \dots +
\alpha_kg_k = \alpha_1g_1u + \dots + \alpha_kg_ku$ in the
alternative loop algebra $F[Q]$. But $F[Q] = FQ/I$. Thus
$\alpha_1g_1 + \dots + \alpha_kg_k = \alpha_1g_1u + \dots +
\alpha_kg_ku (\text{mod} I)$ in the loop algebra $FQ$. We suppose
that the loop $Q$ is infinite. Then, there exist $u \in Q$ such
that $\alpha_1g_1u \notin \{\alpha_1g_1, \ldots, \alpha_kg_k\}$.
By the definition the loop algebra $FQ$ is a free $F$-module with
the basis $\{g \vert g \in Q\}$. Then $\alpha_1g_1u \in I$. From
here it follows that $I = FQ$. But this contradicts Theorem 6.5.
Hence the loop $Q$ is finite.

By \cite[pag. 92]{Bruck}, in the finite Moufang loop $Q$ the order
of any of its element divides the order of $Q$. If $e \neq g \in
Q$ is an element of simple order $p$ then, by item 6) of
Proposition 3.3,  $a = p - (e + g + g^2 + \ldots + g^{p - 1}) \in
\omega[Q]$. We have $(e + g + g^2 + \ldots + g^{p - 1})(e + g +
g^2 + \ldots + g^{p - 1}) = p(e + g + g^2 + \ldots + g^{p - 1})$.
Then $a^2 = p^2 - p(e + g + g^2 + \ldots + g^{p - 1})$ and, by
induction of $n$, it is easy to show that $a^{2^n} = p^{2^n} -
p^{2^n - 1}(e + g + g^2 + \ldots + g^{p - 1})$. We choose an $n$
such that $(\omega[Q])^{2^n} = (0)$. Then $a^{2^n} = 0$. We
suppose that $F$ does not have the characteristic $p$. It follows,
from the equalities $0 = p^{2^n} - p^{2^n - 1}(e + g + g^2 +
\ldots + g^{p - 1}) = p^{2^4 - 1}(p - (e + g + g^2 + \ldots + g^{p
- 1})$, that $0 = p - (e + g + g^2 + \ldots + g^{p - 1})$, $p = e
+ g + g^2 + \ldots + g^{p - 1}$, i.e. $pg = p$, $g = e$. We have
obtained a contradiction as $g \neq e$. Consequently, $F$ has the
characteristic $p$ and $Q$ is $p$-loop. Consequently, 1)
$\Rightarrow$ 2).

 Conversely, let the field $F$ have a
characteristic $p$ and let $Q$ be a  finite $p$-loop. By
Proposition 7.6 it will be centrally nilpotent loop. Let $H = <a>$
be a cyclic group of order $p$ from the center $Z(Q)$ of $Q$. We
will prove that the product
$$(e - a^{i^1})(e - a^{i^2}) \ldots (e - a^{i^m})$$
equals zero, if $m \geq p$. Indeed, if we use the identity $e - xy
= (e - x) + (e - y) - (e - x)(e - y)$, then the last product is
the sum of the factors of type $(e - a)^k$, $k \geq p$. Then
$$(e - a)^p = e - C_p^1a + C_p^2a^2 - \ldots \pm a^p.$$
As all binomial coefficients $C_p^i$ divide by $p$, then they are
zero in the field $F$. Consequently, $(\omega[H])^p = (0)$, where
$\omega[H]$ means the augmentation ideal of the alternative loop
algebra $F[H]$. Let $\mu H$ means the ideal of the alternative
loop algebra $F[Q]$, generated by the set $\{e - h \vert h \in
H\}$. As the subloop $H$ belongs to the center $Z(Q)$, then the
equality $(\omega[H])^p = (0)$ entails the equality $(\mu H)^p =
(0)$.

We will prove the nilpotency of augmentation ideal $\omega[Q]$ via
induction on the order of the loop $Q$. As $H \subseteq Z(Q)$,
then $H$ is normal in $Q$ and $H$ induces the homomorphism $F[Q]
\rightarrow F[Q/H]$. By item 4) of Proposition 3.4, we have
$\omega[Q]/\mu H \cong \omega[Q/H]$. By inductive hypotheses, the
augmentation ideal $\omega[Q/H]$ is nilpotent, for example, of
index $k$. Then $(\omega[Q])^k \subseteq \mu H$ and
$(\omega[Q])^{kp} \subseteq (\mu H)^p = (0)$. Consequently, the
ideal $\omega[Q]$ is nilpotent. Hence 2) $\Rightarrow$ 1).
Consequently,  1) $\Leftrightarrow$ 2).

We use the equivalence of items 1), 2) and  we suppose that 2)
holds. Let $g_1, g_2, \ldots, g_k$ be all elements of $Q$. By item
4 of Proposition 3.3 $\omega Q$ is a finite sum of modules $Fu_i$,
where $u_i = e - g_i$. The field $F$ has a characteristic $p$ and
$g^{p^n} = e$ for some $n$. Then $u_i^{p^n} = (e - g_i)^{p^n} =
0$. Hence $Fu_i$ satisfies the minimum condition for submodules.
It easily follows from here that $\omega Q$ is Artinian, i.e the
item 3) holds. Furthermore, it is known \cite{ZSSS} that the
Zhevlakov radical of an Artinian alternative algebra is nilpotent.
By item 7) of Proposition 3.3, $J(F[Q]) = \omega[Q]$. Thus from 3)
it follows 1). This completes the proof of Theorem 7.7.

\smallskip
Tiraspol State University of Moldova

E-mail: sandumn@yahoo.com

\begin{thebibliography}{99}

\bibitem{Bad} Badeev A. V. \textit{On the Specht property of
varieties of commutative alternative alebras over a field of
characteristic 3 and of commutative Moufang loops,} Sibirsk. Mat.
Zb., 41(2000), 1252 -- 1268 (Russian).

\bibitem{Bruck} Bruck R. H. \textit{A survey of binary systems,}
Berlin-Gottingen-Heidelberg, Springer-Verlag, 1958.

\bibitem{CovSan} Covalschi A. V., Sandu N. I. \textit{On the generalized nilpotent and
generalized solvable loops I,} ROMAI Jurnal, 7, 1(2011), 39 -- 63.

\bibitem{Chein} Chein O. \textit{Moufang loops of small order,}
Trans. Amer. Math. Soc., 188(1974),  31 -- 51.

\bibitem{CPS} Chein O., Pflugfelder H. O., Smith J. D. H. \textit{Quasigroups and
Loops: Theory and applications,} Berlin, Helderman Verlag, 1990.

\bibitem{Dnies} \textit{The Dniester Notebook: Unsolved problems
in the theory of and modules,} Third edition; Akad. Nauk SSSR
Sibirsk Otdel., Inst. Mat., Novosibirsk 1982 (Russian).

\bibitem{Doro} Doro S. \textit{Simple Moufang loops,} Math. Proc. Camb. Phil.
Soc., 83(1978), 377 -- 392.

\bibitem{Fost} Foster A. \textit{The indempotent elements of a
commutative ring form a Boolean algebra,} Duke Math. J., 12(1945),
 143 -- 152.

\bibitem{Glaub} Glauberman G. \textit{On loops of odd order II,}
J. Algebra, 8(1968), 393 -- 414.

\bibitem{Glaub68} Glauberman G., Wright C. R. R. \textit{Nilpotence of
finite Moufang $2$-loops,} J. Algebra, 8(1968), 415 -- 417.




\bibitem{G} Goodaire E. C. \textit{A brief history of loop rings.} 15th Brasilian
School of Algebra (Canela, 1998,) Mat.Contemp., 16(1999), 93--109.

\bibitem{Good87} Goodaire E. C. \textit{Circle loops of radical
alternative rings,} Algebras, groups and geometries, 4(1987), 461
-- 474.


\bibitem{GrShest09} Grishkov A. N, Shestakov I. P. \textit{Commutative
Moufang loops and alternative algebras,} J. Algebra, 333, 1(2011),
1 -- 13; arXiv:0811.3787v1.

\bibitem{GrZ} Grishkov A. N, Zavarnitsine A. V. \textit{Sylow's
theorem for Moufang loops,} J. Algebra, 321, 7, 1(2009), 1813 --
1825; arXiv:0709.2696v1.

\bibitem{Jacob} Jacobson N. \textit{Sructure of rings,} Amer. Math.
Soc. Colloq. Publ., XXXVII, 1956.

\bibitem{Lieb}  Liebeck M. W. \textit{The
classification of finite Moufang loops,} Math. Proc. Camb. Phil.
Soc., 102(1987),  33 -- 47.



\bibitem{Malb} Malbos J.-P. \textit{Sur la classe de nilpotence
des boucles commutative de Moufang et des espaces mediaus,} C.R.
Acad. Sci. Paris Ser. A, 287(1980), 691 -- 693.

\bibitem{Malc} Mal'cev A. I. \textit{About the decomposition of
algebra into a direct sum of radical and semisimple subalgebra,}
Dokl. AN SSSR, 36(1942), 46 -- 50 (Russian).

\bibitem{Manin} Manin Yu. I. \textit{Cubic forms,} Amsterdam:
North-Holland, 1979.

\bibitem{Med} Medvedev Yu. A. \textit{Example of a variety of
alternative algebras over a field of characteristic two,} Algebra
i Logika, 19(1980), 300 -- 313 (Russian).

\bibitem{Mouf} Moufang R.  \textit{Zur Structur von
Alternativek\"orpern,} Math. Ann., 110(1935), 416 -- 430.

\bibitem{Paige} Paige L. J.  \textit{A class of simple Moufang loops,}
Proc. Amer. Math. Soc., 7(1956),  471 -- 482.

\bibitem{Pchel} Pchelintsev S. V. \textit{Structure of finitely
generated commutative alternative algebras and special Moufang
loops,} Matem. zametki, 80(2006), 413 -- 420 (Russian).

\bibitem{Perl} Perlis S. \textit{A characterization of the radical of
an algebra,} Bull. Amer. Math. Soc., 48(1942), 128 -- 132.

\bibitem{Ryab} Ryabukhin Iu. M. \textit{Quasi-regular algebras,
modules, groups and varieties,} Bul. Acad. \c{S}tiin\c{t}e  Repub.
Mold. Mat., 1(1997), 6 -- 62 (Russian).

\bibitem{San87} Sandu N., \textit{Infinite irreducible systems of
identities of commutative Moufang loops and distributive Steiner's
quasigroups,} Izv. Ak. Nauk (USSR), ser. mat., 51, 6(1987), 171 --
189 (Russian).

\bibitem{San99} Sandu N. I. \textit{On the Bruck-Slaby Theorem for
Commutative Moufang Loops,} Matem. zametki, 66, 2(1999), 275 --
281(Russian).


\bibitem{San04} Sandu N. I. \textit{Infinite independent systems of identities
of alternative commutative algebra over a field of characteristic
three,} Discussiones Mathematicae. General Algebra and
Applications, 24(2004), 5 -- 30.

\bibitem{San08} Sandu N. I. \textit{Simple Moufang loops and alternative
algebras,} http://arxiv;1102.1367v1[math.RA].

\bibitem{12} Sandu  N. I. \textit{The classification on sinple
Moufang loops,} http://arxiv:0804.2048v1 [math.GR].

\bibitem{San88} Sandu  N. I. \textit{About the embedded of Moufang
loops in alternative algebras,} http://arxiv:0804.0597 [math. GR].

\bibitem{San99} Sandu N. I.  \textit{Free Moufang loops and alternative
algebras,} Bul. Acad. \c{S}tiin\c{t}e  Repub. Mold. Mat. 3(2009),
 96 - 108.

\bibitem{San11} Sandu N. I. \textit{Simple Moufang loops and Galois extensions,}
http://arxiv;1102.1373v1[math.RA]

\bibitem{Shest04}  Shestakov I. P. \textit{Moufang loops and
alternative algebras,} Proc. Amer. Math. Soc., 132(2004), 313 --
316.


\bibitem{Umir} Umirbaev U. U. \textit{The Specht property of a
variety of solvable alternative algebras,} Algebra i Logika,
24(1985), 226 -- 239 (Russian).

\bibitem{ZSSS} Zhevlakov K.A.,  Slin'ko A. M.,
Shestakov I. P,  Shirshov A. I. \textit{ Rings that are nearly
associative.} Nauka, Moscow, 1978 (Russian); English transl.,
Academic Press, 1982.

\end{thebibliography}
\end{document}